\documentclass{elsart}

\usepackage{amsmath,amssymb,latexsym,epsfig,stmaryrd,tipa}
\usepackage[OT1,T1]{fontenc}

\def\End{\mathop{\rm End}\nolimits}

\def\ker{\mathop{\rm ker}\nolimits}

\def\im{\mathop{\rm Im}\nolimits}

\def\cliff{\mbox{\sl Cliff}}

\def\mod{\mathop{\rm mod}\nolimits}

\def\im{\mathop{\rm Im}\nolimits}

\newcommand{\enorm}[1]{\parallel{#1}\parallel}

\newcommand{\R}{\mathbb{R}}
\newcommand{\C}{\mathbb{C}}
\renewcommand{\H}{\mathbb{H}}
\newcommand{\Z}{\mathbb{Z}}

\renewcommand{\O}{\mathbb{O}}

\newcommand{\mf}{\mathfrak}
\newcommand{\mb}{\mathbf}

\newcommand{\cf}[1]{\hspace{-2pt}{\rm \cite{#1}}\hspace{2pt}}
\newcommand{\db}{\mbox{\textcrd}}

\newcommand{\Db}{\mbox{\DH}}
\newcommand{\Dbb}{\mbox{\textbf{\DH}}}
\renewcommand{\mspace}[1]{\mbox{\hspace{#1}}}
\begin{document}

\thispagestyle{empty}

\begin{center}{\LARGE{\bf
Special metrics and triality}}

\vspace{35pt}

{{\bf Frederik Witt\footnote{present address: NWF I - Mathematik Universit\"at Regensburg, D--93040 Regensburg, F.R.G.}}\\Centre de Math\'ematiques Ecole Polytechnique\\F--91128 Palaiseau, France}

fwitt@math.polytechnique.fr

\vspace{40pt}

{\bf ABSTRACT}

\end{center}

We investigate a new 8--dimensional Riemannian geometry defined by a generic closed and coclosed 3--form with stabiliser $PSU(3)$, and which arises as a critical point of Hitchin's variational principle. We give a Riemannian characterisation of this structure in terms of invariant spinor--valued $1$--forms, which are harmonic with respect to the twisted Dirac operator $\Db$ on $\Delta\otimes\Lambda^1$. We establish various obstructions to the existence of topological reductions to $PSU(3)$. For compact manifolds, we also give sufficient conditions for topological $PSU(3)$--structures that can be lifted to topological $SU(3)$--structures. We also construct the first known compact example of an integrable non--symmetric $PSU(3)$--structure. In the same vein, we give a new Riemannian characterisation for topological quaternionic K\"ahler structures which are defined by an $Sp(1)\cdot Sp(2)$--invariant self--dual 4--form. Again, we show that this form is closed if and only if the corresponding spinor--valued $1$--form is harmonic for $\Db$ and that these equivalent conditions produce constraints on the Ricci tensor.

\bigskip

\newpage
%
%
%
%
%
\section{Introduction}
In the classical setting of Hitchin's variational principle~\cite{hi01}, two structures appear giving rise to a Riemannian metric. In both cases, a critical point is defined by a $3$--form $\rho$ which is harmonic with respect to the metric it induces, i.e.
\begin{equation}\label{harmonic}
d\rho=0,\quad d\star_{\rho}\rho=0.
\end{equation}
In dimension $7$, $\rho$ induces a topological reduction of the frame bundle to a $G_2$--principal fibre bundle, and~(\ref{harmonic}) forces the holonomy of the metric to be contained in $G_2$. In dimension $8$, we get a new type of geometry associated with the structure group $PSU(3)=SU(3)/\ker Ad\subset SO(8)$ which, apart from the initial study in~\cite{hi01}, has been largely unexplored so far.

The most basic example of an $8$--manifold with a $PSU(3)$--structure is $SU(3)$ with $\rho(X,Y,Z)=-B([X,Y],Z)$ built out of the $Ad$--invariant Killing metric $B$. Here, $SU(3)$ is a Riemannian symmetric space and $\rho$ is parallel with respect to the induced Levi--Civita connection. The first problem we tackle is to answer Hitchin's question~\cite{hi01},~\cite{hi04} whether there exist any compact, non--symmetric harmonic $PSU(3)$--manifolds (that is where~(\ref{harmonic}) holds). In this case, the tangent bundle must be associated with a principal $PSU(3)$--bundle with $PSU(3)$ acting in its adjoint representation, and we derive necessary conditions for such a reduction. Since the inclusion $PSU(3)\subset SO(8)$ lifts to $Spin(8)$, any $PSU(3)$--bundle induces a canonical spin structure, so the underlying manifold is spin. More importantly, half of the tangent bundle trivialises, that is, there exist four pointwise linearly independent vector fields. It follows that any compact homogeneous Riemannian manifold $G/H$ with $G$ simple is diffeomorphic to $SU(3)$ (Proposition~\ref{topclass}). The question of finding sufficient conditions for $PSU(3)$--bundles to exist on connected, compact spin manifolds is, unlike the analogous question for $G_2$, rather involved. To give at least a partial answer, we restrict ourselves to $PSU(3)$--bundles with vanishing {\em triality class} (Theorem~\ref{psu3structure}). This class is the cohomological obstruction for lifting the structure group from $PSU(3)$ to $SU(3)$, which can be thought of as the analogue of lifting an orthonormal frame bundle to a spin structure. Over 4--dimensional manifolds, this issue has been addressed in~\cite{wo82b} motivated by the fact that the group $SU(3)$ acting in its adjoint representation is the gauge group of quantum chromodynamics~\cite{gene64}. As a result, we are left with rather severe restrictions on the topology of the underlying manifold. In fact, all our examples of harmonic $PSU(3)$--structures have trivial tangent bundle: We build a compact non--symmetric example out of a special $6$--dimensional nilmanifold times a $2$--torus. We also find a family of {\em local} examples obtained from a $4$--dimensional hyperk\"ahler manifold times flat Euclidean $4$--space.

Motivated by the $G_2$--case, we also approach $PSU(3)$--manifolds from a Riemannian point of view and ask: What is the extra datum reducing the orthonormal frame bundle of a Riemannian $8$--manifold to a $PSU(3)$--bundle? For $G_2$, this is a nonvanishing spinor field, and the $G_2$--structure is harmonic if and only if the corresponding spinor field is parallel with respect to the Levi--Civita connection. In~\cite{hi01} Hitchin remarks that, for a $PSU(3)$--structure, there exist two invariant spinor--valued $1$--forms $\sigma_{\pm}\in\Delta_{\pm}\otimes\Lambda^1$. He then shows (albeit with some minor mistakes, cf. Remark~\ref{proof}) that under~(\ref{harmonic}), these are harmonic with respect to the twisted Dirac operator $\Db_{\pm}:\Gamma(\Delta_{\pm}\otimes\Lambda^1)\to\Gamma(\Delta_{\mp}\otimes\Lambda^1)$, that is, $\Db_{\pm}(\sigma_{\pm})=0$. We prove the converse -- this is where triality comes in. The vector representation $\Lambda^1$ and the two irreducible spin representations $\Delta_+$ and $\Delta_-$ of $Spin(8)$ are, though inequivalent as $Spin(8)$--modules, isomorphic as Euclidean vector spaces. It is therefore sufficient to consider the set of isometries $\sigma:\Delta_-\to\Delta_+$ such that $\sigma$ lies in an irreducible subspace of $\Delta_+\otimes\Delta_-$. After working out the $Spin(8)$--orbit structure on this set (Theorem~\ref{orbitdecomposition}), we see that one orbit is isomorphic to $Spin(8)/\big(PSU(3)\times\Z_2\big)$ and corresponds to orientation--{\em preserving} isometries $\sigma_{\pm}:\Lambda^1\to\Delta_{\pm}$ in $\ker\mu_{\pm}\subset\Delta_{\pm}\otimes\Lambda^1$, the kernel of Clifford multiplication. Further, harmonicity with respect to $\Db$ enforces~(\ref{harmonic}).

Surprisingly, another orbit of interest shows up: $Spin(8)/Sp(1)\cdot Sp(2)$, where $Sp(1)\cdot Sp(2)$ stabilises an orientation--{\em reversing} isometry $\sigma_+:\Lambda^1\to\Delta_+$ in $\ker\mu_+$. This yields a new Riemannian characterisation of $Sp(1)\cdot Sp(2)$--structures which so far in the literature have been defined in terms of an invariant self--dual $4$--form $\Omega$, following the higher--dimensional analogy with $Sp(1)\cdot Sp(k)$--structures on $M^{4k}$. A Riemannian manifold $M^{4k}$ whose holonomy is contained in $Sp(1)\cdot Sp(k)$ is called {\em quaternionic K\"ahler}~\cite{sa82} and is necessarily Einstein. It is known that for $k\geq3$, this is equivalent to $d\Omega=0$~\cite{sw89}, that is,~(\ref{harmonic}) holds for $\rho=\Omega$. However, there are counterexamples for $k=2$~\cite{sa01}. Here,~(\ref{harmonic}) holds if and only if $\Db_+(\sigma_+)=0$. Although this does not imply that the metric is Einstein, as we will show by using Salamon's counterexample, we nevertheless obtain non--trivial constraints on the Ricci tensor (Proposition~\ref{riccitensor}).

Finally, we remark on the relationship with other distinguished Riemannian metrics. A third characterisation of $PSU(3)$--structures is given in terms of a symmetric $3$--tensor. This fits into a series of special ``nearly--integrable'' Riemannian metrics which were investigated in~\cite{nu06}. Although this integrability condition is in a way opposite to ours (cf. Theorem~\ref{psu3thm} and Remark~\ref{nearlyint}), it links $PSU(3)$--structures to $SO(3)$--structures in dimension $5$~\cite{bonu07}, matching the relationship between $G_2$--manifolds in dimension $7$ and hyperk\"ahler manifolds in dimension $4$. This does not only indicate a way to construct further examples of harmonic $PSU(3)$--structures, but also provides evidence for a still unexplored, intrinsic relationship between these special Riemannian geometries in low dimensions.
%
%
%
%
%
\section{Triality and supersymmetric maps}
\label{triality}
In the presence of a metric, we can identify vectors in $\R^8$ with $1$--forms in $\Lambda^1=\Lambda^1\R^{8*}$ and we shall freely do so throughout this paper. The triality principle asserts that the vector representation $\pi_0:Spin(8)\to SO(\Lambda^1)$ and the two chiral spin representations
$\pi_{\pm}:Spin(8)\to SO(\Delta_{\pm})$ are isomorphic as Euclidean vector spaces even though they are inequivalent as irreducible $Spin(8)$--spaces. More precisely, the representations are related by $\pi_0=\pi_+\circ\kappa\circ\lambda$ and $\pi_-=\pi_+\circ\lambda^2$, where $\kappa$ and $\lambda$ are two outer $Spin(8)$--automorphisms of order two and three. Morally, this means that we can exchange any two of the representations $\Lambda^1$, $\Delta_+$ and $\Delta_-$ by an outer automorphism, while the remaining third one is fixed. A convenient model for the underlying Euclidean vector space is provided by the octonions $\O$. Here, $Spin(8)$ acts as orientation preserving isometry group of the inner product induced by the oriented orthonormal basis 
\begin{equation}\label{octbasis}
1,i,j,k,e,e\cdot i,e\cdot j,e\cdot k.
\end{equation} 
If $R_u$ denotes right multiplication by $u\in\O$, the map
\begin{equation}\label{cliff(8)}
u\in\O\mapsto \left(\begin{array}{cc} 0 & R_u\\ -R_{\bar{u}} & 0
\end{array}\right)\in\End(\O\oplus\O)
\end{equation}
extends to an isomorphism $\cliff(\O)\cong\End(\O\oplus\O)$ where $\Delta=\O\oplus\O$ is the (reducible) space of spinors for $Spin(8)$. These two summands can be distinguished by an orientation, since a volume form acts on these by $\pm Id$, which gives rise to the spin representations $\Delta_+$ and $\Delta_-$. The explicit matrix representation~(\ref{cliff(8)}) we will use throughout this paper is given in Appendix~\ref{matrix}. Moreover, the inner product on $\O$ can be adopted as the $Spin(8)$--invariant inner product $q$ on $\Delta_{\pm}$.

\begin{defn}
A {\em supersymmetric map} is an isometry between two of the three spaces $\Lambda^1$, $\Delta_+$ or $\Delta_-$, which lies in an irreducible $Spin(8)$--submodule of $\Lambda^1\otimes\Delta_{\pm}$, $\Delta_{\pm}\otimes\Lambda^1$ or $\Delta_{\pm}\otimes\Delta_{\mp}$.
\end{defn}

\begin{exmp}
A unit spinor $\Psi\in\Delta_+$ induces a supersymmetric map $X\in\Lambda^1\mapsto X\cdot\Psi\in\Delta_-$. As an element in $\Delta_-\otimes\Lambda^1\cong\Delta_+\oplus\ker\mu_-$, where $\mu_{\pm}:\Delta_{\pm}\otimes\Lambda^1\to\Delta_{\mp}$ denotes Clifford multiplication, it belongs to the irreducible subspace $\Delta_+$. One easily checks that $Spin(8)$ acts transitively on the set of supersymmetric maps in $\Delta_+$, and that the orbit is isomorphic with $Spin(8)/Spin(7)$. In passing we remark that $\ker\mu_{\pm}\cong\Lambda^3\Delta_{\mp}$.
\end{exmp}

The jargon has its origin in particle physics where a supersymmetry is supposed to transform bosons (particles which are elements in a vector representation of the spin group) into fermions (particles which are elements in a spin representation of the spin group).

The case of supersymmetric maps which are induced by a 3--form over $\Lambda^1$, $\Delta_+$ or $\Delta_-$ is more interesting, and we set out to give a complete classification. As we are only concerned with the metric structure of these spaces, triality implies that we are free to consider the module $\Delta_+\otimes\Delta_-$ rather than $\Delta_{\pm}\otimes\Lambda^1$, and we subsequently do so for various reasons. As a $Spin(8)$--module, $\Delta_+\otimes\Delta_-\cong\Lambda^1\oplus\Lambda^3$, and we define
$$
\mf{I}_g=\{\rho\in\Lambda^3\subset\Delta_+\otimes\Delta_-\,|\,\rho:\Delta_-\to\Delta_+\mbox{ is a supersymmetric map}\}.
$$
This set is acted on by $Spin(8)$ and we exhibit the orbit structure based on the following

\begin{thm}\label{orbit}
If $\rho\in\Lambda^3$ lies in $\mf{I}_g$, then $\rho$ is of unit length
and there exists a Lie bracket $[\cdot\,,\cdot]$ on $\Lambda^1$ such
that
\begin{equation}\label{bracket}
\rho(x,y,z)=g([x,y],z).
\end{equation}
Consequently, the adjoint group of this Lie algebra acts as a
group of isometries on $\Lambda^1$.

Conversely, if there exists a Lie algebra structure on $\Lambda^1$
whose adjoint group leaves $g$ invariant, the 3--form defined by
(\ref{bracket}) and divided by its norm belongs to $\mf{I}_g$.
\end{thm}

\begin{pf}
Because of the skew--symmetry of $\rho$, the metric $g$ is necessarily invariant
under the adjoint action of the induced Lie algebra, for
$$
g([x,y],z)=\rho(x,y,z)=-g([x,z],y).
$$
Being an isometry inducing a Lie bracket through (\ref{bracket}) and vice versa are both quadratic conditions on the coefficients of $\rho$ which we show to coincide. We define the linear map
$$
Jac:\Lambda^3\otimes\Lambda^3\to\Lambda^4
$$
by skew--symmetrising the contraction to $\Lambda^2\otimes\Lambda^2$.
This is most suitably expressed in index notation with respect to some orthonormal basis $\{e_i\}$, namely
\begin{eqnarray*}
Jac(\rho_{ijk}\tau_{lmn}) & = &
\rho_{[ij}^{\phantom{[ij}k}\tau_{lm]k}\\
& = &
\frac{1}{6}\big(\rho_{ij}^{\phantom{ij}k}\tau_{lmk}+\rho_{il}^{\phantom{il}k}\tau_{mjk}+\rho_{im}^{\phantom{im}k}\tau_{jlk}+\rho_{jl}^{\phantom{jl}k}\tau_{imk}+\rho_{jm}^{\phantom{jm}k}\tau_{lik}\\
& & +\rho_{lm}^{\phantom{lm}k}\tau_{ijk}\big).
\end{eqnarray*}
In particular,
\begin{equation}\label{jacid}
Jac(\rho_{ijk}\rho_{lmn})=\frac{1}{3}\big(\rho_{ij}^{\phantom{ij}k}\rho_{klm}+\rho_{li}^{\phantom{li}^k}\rho_{kjm}+\rho_{jl}^{\phantom{jl}k}\rho_{kim}\big).
\end{equation}
If we are given a 3--form $\rho$ and define a skew--symmetric map
$[\cdot\,,\cdot]:\Lambda^2\to\Lambda^1$ by (\ref{bracket}), then
the Jacobi identity holds, i.e. we have defined a Lie bracket, if and only if $Jac(\rho\otimes\rho)=0$.

Next we analyse the conditions for $\rho$ to induce an isometry. For a $p$--form $\rho$ we have $q(\rho\cdot\Psi_1,\Psi_2)=(-1)^{p(p+1)/2}q(\Psi_1,\rho\cdot\Psi_2)$, so $\rho$ defines an isometry $\Delta_{\pm}\to\Delta_{\mp}$ if and only if for any pair of spinors of equal chirality, $q(\rho\cdot\rho\cdot\Psi_1,\Psi_2)=q(\Psi_1,\Psi_2)$ holds. Considering the $Spin(8)$--equivariant maps
$$
\Gamma_{\pm}:\rho\otimes\tau\in\Lambda^3\otimes\Lambda^3\mapsto\rho\cdot\tau\in\cliff(\Lambda^1)\cong\End(\Delta)\stackrel{pr_{\pm}}{\mapsto}\rho\cdot\tau_{\,|\Delta_{\pm}}\in\Delta_{\pm}\otimes\Delta_{\pm},
$$
this condition reads $\rho\in\mf{I}_g$ if and only
if $\Gamma_{\pm}(\rho\otimes\rho)={\rm Id}_{\Delta_{\pm}}$. Using
the algorithm in~\cite{sa89} or a suitable computer programme, we decompose both the domain and the target space into irreducible components,
\begin{eqnarray*}
\Lambda^3\otimes\Lambda^3 & \cong &
\mathbf{1}\oplus2\Lambda^2\oplus\Lambda^4_+\oplus\Lambda^4_-\oplus
[1,4,3,3]\oplus[2,4,2,3]\oplus[2,4,3,2]\oplus\\
& & [2,4,2,2]\oplus2[2,3,2,2]\oplus[2,2,1,1]\\
\Delta_{\pm}\otimes\Delta_{\pm}& = & \Lambda^2\Delta_{\pm}\oplus\odot^2\Delta_{\pm}\cong\Lambda^2\oplus\mathbf{1}\oplus\Lambda^4_{\pm},
\end{eqnarray*}
where we label irreducible representations by their highest weight (expressed in the basis of fundamental roots). The modules $\Lambda^4_+=[1,2,2,1]$ and $\Lambda^4_-=[1,2,1,2]$ are
the spaces of self--dual and anti--self--dual 4--forms respectively.
Note that $\Gamma_+(\rho\otimes\tau)^{tr}=\Gamma_-(\tau\otimes\rho)$
and so it suffices to consider the map $\Gamma_+$ only. Since the
map induced by $\rho$ is symmetric, it follows that
$\Gamma_{\pm}(\rho\otimes\rho)\in\odot^2\Delta_{\pm}=\mathbf{1}\oplus\Lambda^4_{\pm}$.
Moreover the image clearly contains $\Lambda^4_{\pm}$. As a result, $\Gamma_{\pm}(\rho\otimes\rho)_{\bigodot^2_0\Delta_{\pm}}=0$ is
a necessary condition for $\rho$ to lie in $\mf{I}_g$.

Next we identify this obstruction in $\Lambda^4_+$ with $
Jac(\rho\otimes\rho)$ by showing
\begin{equation}\label{isom++}
\Gamma_+(\rho\otimes\rho)\oplus\Gamma_-(\rho\otimes\rho)=-3
Jac(\rho\otimes\rho)_{\Lambda^4_+}+\enorm{\rho}^2\mathop{\rm Id}.
\end{equation}
We first remark that Clifford multiplication induces a map
$$
\rho\otimes\tau\in\Lambda^3\otimes\Lambda^3\mapsto\rho\cdot\tau\in\Lambda^0\oplus\Lambda^2\oplus\Lambda^4\oplus\Lambda^6
$$
if we regard the product $\rho\cdot\tau$ as an element of
$\cliff(\Lambda^1,g)\cong\Lambda^*$ under the natural isomorphism. The various
components of $\rho\cdot\tau$ under this identification are accounted for by the ``coinciding pairs`` (c.p.) in the expression
$\rho_{ijk}\tau_{lmn}e_{ijklmn}$, $i<j<l$, $l<m<n$. For instance, having three coinciding pairs implies
$i=l$, $j=m$ and $k=n$, hence $e_{ijklmn}=1$. Then
$\rho=\sum\nolimits_{i<j<k}c_{ijk}e_{ijk}$ gets mapped to
\begin{eqnarray*}
\rho\cdot\rho & = & \sum_{\begin{array}{c}\scriptstyle
i<j<k\\\scriptstyle
l<m<n\end{array}}\rho_{ijk}\rho_{lmn}e_{lmnijk}\\
& = & \sum_{\begin{array}{c}
\scriptstyle i<j<k\\\scriptstyle l<m<n\end{array}\mbox{\scriptsize ,
3 c.p.}}\rho_{ijk}\rho_{lmn}e_{lmnijk} +\sum_{\begin{array}{c}
\scriptstyle i<j<k\\\scriptstyle l<m<n\end{array}\mbox{\scriptsize ,
1 c.p.}}\rho_{ijk}\rho_{lmn}e_{lmnijk}.
\end{eqnarray*}
There is no contribution by the sum of two c.p. as $\rho\cdot\rho$ is symmetric. Now the first sum is just
$$
\enorm{\rho}^2\cdot\mathbf{1}=\sum\nolimits_{i<j<k}\rho^2_{ijk}\mathbf{1}
$$
which leaves us with the contribution of the sum with one pair of equal indices. No matter which indices of the two triples $(i<j<k)$ and $(l<m<n)$ coincide, the skew--symmetry of the
$c_{ijk}$ and $e_{ijk}$ allows us to rearrange and rename the indices in such a way that the second sum equals
\begin{eqnarray*}
\sum_a\sum_{\begin{array}{c} \scriptstyle
j<k,\,m<n\\\mbox{\scriptsize{\it j,k,m,n}
dist.}\end{array}}\rho_{ajk}\rho_{amn}e_{amnajk} & = &
-\sum_a\sum_{\begin{array}{c} \scriptstyle
j<k,\,m<n\\\mbox{\scriptsize{\it j,k,m,n}
dist.}\end{array}}\rho_{ajk}\rho_{amn}e_{mnjk}\\
& = & -3Jac(\rho\otimes\rho).
\end{eqnarray*}
This implies~(\ref{isom++}), hence the assertion of the theorem.
\end{pf}

Consequently, the 3--forms in $\mf{I}_g$ encode a Lie algebra structure whose adjoint action preserves the metric on $\Lambda^1$. We also say that the Lie structure is {\em adapted} to the metric $g$ and write $\mf{l}$ if we think of $\Lambda^1$ as a Lie
algebra. We classify the resulting Lie algebras next.

Let us recall some basic notions (see for instance~\cite{onvi94}). A Lie algebra $\mf{g}$ is said to be {\em simple} if it contains no non--trivial ideals. A {\em semi--simple} Lie algebra is a direct sum of simple ones which is to say that it does not possess any non--trivial abelian ideal.  Equivalently, $\mf{g}^{(1)}=[\mf{g},\mf{g}]=\mf{g}$.
On the other hand, if the derived series defined inductively by $\mf{g}^{(k)}=[\mf{g}^{(k-1)},\mf{g}^{(k-1)}]$ becomes trivial from some integer $k$ on, then $\mf{g}$ is {\em solvable}. Any abelian Lie algebra is solvable and so is any sub--algebra of a solvable one. Moreover, every Lie algebra contains a maximal solvable ideal, the so--called {\em radical} $\mf{r}(\mf{g})$ of $\mf{g}$. In particular, the centre $\mf{z}(\mf{g})$ is contained in $\mf{r}(\mf{g})$. If there is equality, then $\mf{g}$ is said to be {\em reductive}. Reductive Lie algebras are a direct Lie algebra sum of their centre and a semi--simple Lie algebra.

\begin{prop}\label{reductive}
An adapted Lie algebra $\mf{l}$ is reductive.
\end{prop}

\begin{pf}
By the lemma below, $\mf{r}(\mf{g})$ is abelian which implies $g([R_1,X],R_2)=-g(X,[R_1,R_2])=0$
for any $X\in\mf{l}$ and
$R_1,\,R_2\in\mf{r}(\mf{g})$. Therefore $[X,R_1]\in\mf{r}(\mf{g})\cap\mf{r}(\mf{g})^{\perp}=\{0\}$, hence $\mf{r}(\mf{g})=\mf{z}(\mf{g})$.
\end{pf}

\begin{lem}
Let $\mf{s}$ be a solvable Lie algebra which is adapted to some metric $g$. Then $\mf{s}$ is abelian.
\end{lem}

\begin{pf}
We proceed by induction over $n$, the dimension of $\mf{s}$.
If $n=1$, then $\mf{s}$ is abelian and the assertion is trivial. Now
assume that the assumption holds for all $1\leq m<n$. Let $\mf{a}$
be a non--trivial abelian ideal of $\mf{s}$. This, of course, does exist, for otherwise $\mf{s}$ would be semi--simple. The $ad$--invariance of $g$ implies $g([A,X],Y)=0$ for all $X,\,Y\in\mf{s}$ and $A\in\mf{a}$.
For if $X\in\mf{a}$, then $[A,X]=0$ and if $X\in\mf{a}^{\perp}$, then $g([A,X],Y)=-g(X,[A,Y])=0$
since $[A,Y]\in\mf{a}$. Hence $\mf{a}\subset\mf{z}(\mf{s})$. We can therefore split $\mf{s}=\mf{z}\oplus\mf{h}$ into a direct sum of vector spaces with $\mf{h}$ an
orthogonal complement to $\mf{z}$ of dimension strictly less than $n$. Now for all $X\in\mf{s}$,
$Z\in\mf{z}$ and $H\in\mf{h}$ we have $g([X,H],Z)=-g(H,[X,Z])=0$,
so $[X,H]\in\mf{z}^{\perp}=\mf{h}$, or equivalently, $\mf{h}$ is
an ideal of $\mf{s}$. As such, it is adapted and solvable since $\mf{s}$ is adapted and solvable. Hence the induction hypothesis applies and $\mf{s}$ is abelian.
\end{pf}

As a result, we are left to determine the semi--simple part of an adapted Lie algebra $\mf{l}$ of dimension $8$. Appealing to Cartan's classification of simple Lie algebras, we obtain the following possibilities ($\mf{z}^p$ denoting the centre of dimension $p$):
\begin{enumerate}
\item $\mf{l}_1=\mf{su}(3)$
\item $\mf{l}_2=\mf{su}(2)\oplus\mf{su}(2)\oplus \mf{z}^2$
\item $\mf{l}_3=\mf{su}(2)\oplus\mf{z}^5$.
\end{enumerate}
Hence there is a disjoint decomposition of $\mf{I}_g$ into the sets
$\mf{I}_{g1}$, $\mf{I}_{g2}$ and $\mf{I}_{g3}$ acted on by $Spin(8)$
and pooling together the forms which induce the Lie algebra structure $\mf{l}_1$, $\mf{l}_2$ or $\mf{l}_3$.

\begin{thm}\label{orbitdecomposition}
The sets $\mf{I}_{g1}$, $\mf{I}_{g2}$ and $\mf{I}_{g3}$ can be
described as follows:
\begin{enumerate}
\item $\mf{I}_{g1} = Spin(8)/\big(PSU(3)\times\Z_2\big)$
\item $\mf{I}_{g2} = (0,1)\times Spin(8)/\big(SU(2)\cdot SU(2)\times U(1)\big)$
\item $\mf{I}_{g3} = Spin(8)/Sp(1)\cdot Sp(2)$,
\end{enumerate}
where $SU(2)\cdot SU(2)=SU(2)\times SU(2)/\Z_2$ and $Sp(1)\cdot Sp(2)=Sp(1)\times Sp(2)/\Z_2$ cover the standard inclusions $SO(3)\times SO(3)\hookrightarrow SO(8)$ and $SO(3)\times SO(5)\hookrightarrow SO(8)$. Further, we have the $Spin(8)$--invariant decomposition $\mf{I}_g=\mf{I}_{g+}\cup\mf{I}_{g-}$ into 3--forms whose induced isometry is orientation--preserving or --reversing respectively. Then $\mf{I}_{g-}=Spin(8)/\big(PSU(3)\times\Z_2\big)$
and $\mf{I}_{g+}$ foliates over the circle $S^1$ with principal
orbits $Spin(8)/\big(SU(2)\cdot SU(2)\times U(1)\big)$ over $S^1-\{pt\}$ and
a degenerate orbit $Spin(8)/Sp(1)\cdot Sp(2)$ at $\{pt\}$.
\end{thm}

\begin{pf}
We remark that the stabiliser of $\rho_i\in\mf{I}_{gi}$ in $SO(8)$ is $SO(8)\cap Aut(\mf{l}_i)$. Consider first a 3--form $\rho_1\in\mf{I}_{g1}$, that is, $\rho_1$ induces an $\mf{su}(3)$--structure on $\Lambda^1$. Since the fixed Riemannian metric $g$ is $ad$--invariant it must coincide with the (negative definite) Killing form $B(X,Y)={\rm Tr}(ad_X\circ ad_Y)$ up to a negative constant $c$. It is well known (cf. for instance~\cite{gene64}) that there exists an orthogonal basis $e_1,\ldots,e_8$ such that the totally anti--symmetric structure constants $c_{ijk}$ are given by
$$
c_{123}=1,\quad c_{147}=-c_{156}=c_{246}=c_{257}= c_{345}=-c_{367}=\frac{1}{2},\quad c_{458}=c_{678}=\frac{\sqrt{3}}{2}
$$
and $B(e_i,e_i)=-3$. Hence $f_i=e_i/\sqrt{-3c}$ is $g$--orthonormal. The relation~(\ref{bracket}) and the requirement to be of unit norm implies that
\begin{eqnarray*}
\rho_1 & = & \frac{1}{\sqrt{-3c}}f_{123}+\frac{1}{2\sqrt{-3c}}f_1(f_{47}-f_{56})+\frac{1}{2\sqrt{-3c}}f_2(f_{46}+f_{57})\\
& & +\frac{1}{2\sqrt{-3c}}f_3(f_{45}-f_{67})+\frac{1}{2\sqrt{-c}}f_8(f_{45}+f_{67})\\
& = & \frac{1}{2}f_{123}+\frac{1}{4}f_1(f_{47}-f_{56})+\frac{1}{4}f_2(f_{46}+f_{57})+\frac{1}{4}f_3(f_{45}-f_{67})\\
& &+\frac{\sqrt{3}}{4}f_8(f_{45}+f_{67})
\end{eqnarray*}
where as usual, the notation $f_{ijk}$ will be shorthand for
$f_i\wedge f_j\wedge f_k$ and vectors are identified with their dual
in presence of a metric. Any 3--form of $\mf{I}_{g1}$ being representable in this way, it follows that $SO(8)$ acts transitively on $\mf{I}_{g1}$. The stabiliser of $\rho_1$ in $SO(8)$ is the adjoint group $SU(3)/\Z_3=PSU(3)$. As $\pi_1\big(PSU(3)\big)=\Z_3$, this is covered by $PSU(3)\times\Z_2$ in $Spin(8)$, hence $\mf{I}_{g1}=Spin(8)/\big(PSU(3)\times\Z_2\big)$. Using the matrix representation of $\cliff(\Lambda^1)$ given in Appendix~\ref{matrix} with respect to some ordered basis $\Psi_{i\pm}$ of $\Delta_{\pm}$, the isometry $A_{\rho_1}:\Delta_-\to\Delta_+$ induced by $\rho_1$ is
\begin{equation}\label{rhomap}
A_{\rho_1}=\frac{1}{4}\left
(\begin {array}{rrrrrrrr}
\scriptstyle\sqrt{3}&\scriptstyle0&\scriptstyle0&\scriptstyle3&\scriptstyle-\sqrt {3}&\scriptstyle0&\scriptstyle0&\scriptstyle1\\
\scriptstyle 2 &\scriptstyle-\sqrt{3}&\scriptstyle-1&\scriptstyle0&\scriptstyle2&\scriptstyle-\sqrt{3}&\scriptstyle-1&\scriptstyle0\\
\scriptstyle 0&\scriptstyle3&\scriptstyle-\sqrt{3}&\scriptstyle0&\scriptstyle0&\scriptstyle-1&\scriptstyle-\sqrt {3}&\scriptstyle0\\
\scriptstyle -1&\scriptstyle0&\scriptstyle2&\scriptstyle\sqrt{3}&\scriptstyle1&\scriptstyle0&\scriptstyle-2 &\scriptstyle-\sqrt{3}\\
\scriptstyle-\sqrt{3}&\scriptstyle0&\scriptstyle0&\scriptstyle1&\scriptstyle\sqrt{3}&\scriptstyle0&\scriptstyle0&\scriptstyle3\\
\scriptstyle-2&\scriptstyle-\sqrt{3}&\scriptstyle-1&\scriptstyle0&\scriptstyle-2& \scriptstyle-\sqrt{3}&\scriptstyle-1&\scriptstyle0\\
\scriptstyle0&\scriptstyle-1&\scriptstyle-\sqrt{3}&\scriptstyle0&\scriptstyle0&\scriptstyle3&\scriptstyle-\sqrt{3}&\scriptstyle0\\
\scriptstyle1&\scriptstyle0&\scriptstyle2&\scriptstyle-\sqrt{3}&\scriptstyle-1&\scriptstyle0&\scriptstyle-2&\scriptstyle\sqrt{3}
\end {array}\right),
\end{equation}
hence $\det(A_{\rho_1})=-1$. Moreover, we have $\det\pi_{\pm}(a)=1$
for any $a\in Spin(8)$ as the generators $e_i\cdot e_j$ square to $-Id$
and are therefore of determinant 1. The $Spin(8)$--equivariance of
the embedding $\Lambda^3\to\Delta\otimes\Delta$ entails $A_{\pi_0(a)^*\rho_1}=\pi_+(a)\circ A_{\rho_1}\circ\pi_-(a)^{-1}$,
whence $\mf{I}_{g1}\subset\mf{I}_{g-}$. 

Next we turn to the Lie algebras $\mf{l}_2$ and $\mf{l}_3$ where the latter can be seen as a degeneration of the former. So assume $\rho_2$ to be an element of $\mf{I}_{g2}$ inducing an $\mf{su}(2)\oplus\mf{su}(2)\oplus\mf{z}^2$--structure. The restriction to $g$ to any copy of $\mf{su}(2)$ must be as above a negative multiple of the Killing form of $\mf{su}(2)$, so $g=c_1B_1\oplus c_2B_2\oplus g_{|\mf{z}^p}$. There exists a basis $e_i$ of $\mf{su}(2)$ such that $[e_i,e_j]=\epsilon_{ijk}e_k$ (where $\epsilon_{ijk}$ is totally anti--symmetric) and $B(e_i,e_i)=-2$. Choosing such a basis for each copy of $\mf{su}(2)$ and extending this to an orthonormal basis $f_i$ of $\Lambda^1$ by normalising, the requirement on $\rho_2$ to be of unit norm implies
\begin{eqnarray*}
\rho_2 & = & \frac{1}{\sqrt{-2c_1}}f_{123}+\frac{1}{\sqrt{-2c_1}}f_{456}\\
& = & \frac{1}{\sqrt{-2c_1}}f_{123}+\sqrt{\frac{2c_1+1}{2c_1}}f_{456},
\end{eqnarray*}
where $c_1=-\sin(\pi\alpha)/2$, $\alpha\in(0,1)$ is the only $SO(8)$--invariant of $\rho_2$. It follows that $\mf{I}_{g2}$ foliates in $SO(8)$--orbits over $(0,1)$. The automorphism group is $SU(2)/\Z_2\times SU(2)/\Z_2\times GL(2)=SO(3)\times SO(3)\times GL(2)$ and since the Lie algebra structure is adapted to $g$, the stabiliser of $\rho_2$ in $SO(8)$ is given by $SO(3)\times SO(3)\times SO(2)$. This is covered twice by $SU(2)\cdot SU(2)\times U(1)\subset Spin(8)$ and we obtain $\mf{I}_{g2}=(0,1)\times Spin(8)/\big(SU(2)\cdot SU(2)\times U(1)\big)$. The induced isometry $\Delta_-\to\Delta_+$ is
$$
A_{\rho_2}=\left(\begin {array}{rrrrrrrr} \scriptstyle0&\scriptstyle0&\scriptstyle\sqrt {\frac{2c_1+1}{2c_1}} &\scriptstyle\frac {1}{\sqrt {-2c1}}&\scriptstyle0&\scriptstyle0&\scriptstyle0&\scriptstyle0
\\
\scriptstyle0&\scriptstyle0&\scriptstyle-\frac {1}{\sqrt {-2c1}}&\scriptstyle\sqrt {\frac{2c_1+1}{2c_1}} &\scriptstyle0&\scriptstyle0&\scriptstyle0&\scriptstyle0
\\
\scriptstyle\sqrt {\frac{2c_1+1}{2c_1}}&\scriptstyle\frac {1}{\sqrt {-2c1}}&\scriptstyle0&\scriptstyle0&\scriptstyle0&\scriptstyle0&\scriptstyle0&\scriptstyle0
\\
-\scriptstyle\frac {1}{\sqrt {-2c1}}&\scriptstyle\sqrt {\frac{2c_1+1}{2c_1}}&\scriptstyle0&\scriptstyle0&\scriptstyle0&\scriptstyle0&\scriptstyle0&\scriptstyle0
\\
\scriptstyle0&\scriptstyle0&\scriptstyle0&\scriptstyle0&\scriptstyle0&\scriptstyle0&\scriptstyle\sqrt {\frac{2c_1+1}{2c_1}}&\scriptstyle\frac {1}{\sqrt {-2c1}}
\\
\scriptstyle0&\scriptstyle0&\scriptstyle0&\scriptstyle0&\scriptstyle0&\scriptstyle0&-\scriptstyle\frac {1}{\sqrt {-2c1}}&\scriptstyle\sqrt {\frac{2c_1+1}{2c_1}}
\\
\scriptstyle0&\scriptstyle0&\scriptstyle0&\scriptstyle0&\scriptstyle\sqrt {\frac{2c_1+1}{2c_1}}&\scriptstyle\frac {1}{\sqrt {-2c1}}&\scriptstyle0&\scriptstyle0
\\
\scriptstyle0&\scriptstyle0&\scriptstyle0&\scriptstyle0&-\scriptstyle\frac {1}{\sqrt {-2c1}}&\scriptstyle\sqrt {\frac{2c_1+1}{2c_1}}&\scriptstyle0&\scriptstyle0
\end {array}\right)
$$
and thus of positive determinant. We conclude as above that $\mf{I}_{g2}\subset\mf{I}_{g+}$.

We obtain the last case for $c_2=0$, i.e. $c_1=-1/2$. Here the stabiliser in $SO(8)$ is isomorphic to $SO(3)\times SO(5)$ whose double cover to $Spin(8)$ is $Sp(1)\cdot Sp(2)$ (using the isomorphisms between $SU(2)\cong Sp(1)$ and $Spin(5)\cong Sp(2)$). Moreover, $\mf{I}_{g3}\subset\mf{I}_{g+}$, whence the theorem.
\end{pf}

By the triality principle, we can exchange $\Delta_+$ or $\Delta_-$ with $\Lambda^1$ while leaving $\Delta_-$ or $\Delta_+$ fixed. Hence we get an analogous orbit decomposition for $\Delta_{\pm}\otimes\Lambda^1$ where the stabiliser subgroups sit now in $SO(\Delta_{\mp})$ and lift via $\pi_{\mp}$ to $Spin(8)$. Note however that the  characterisation of $\mf{I}_{g\pm}$ does depend on the module under consideration as the outer triality morphisms reverse the orientation. In any case, the covering group in $Spin(8)$ acts on all three representations and we analyse now this action in detail. Again it suffices to discuss the case where the stabiliser of the isometry lifts via $\pi_0$.

We start with the group $PSU(3)\times \Z_2$ which projects to $PSU(3)$ in $SO(\Lambda^1)$, $SO(\Delta_+)$ and $SO(\Delta_-)$. Hence $PSU(3)\subset SO(\Lambda^1)$ also gives rise to $PSU(3)$--invariant isometries in $\Delta_{\pm}\otimes\Lambda^1$. We immediately deduce that restricted to $PSU(3)$ in $Spin(8)$, the representation spaces $\Lambda^1$, $\Delta_+$ and
$\Delta_-$ are equivalent. In particular, Clifford multiplication
$\mu_{\pm}:\Lambda^1\otimes\Delta_{\pm}\to\Delta_{\mp}$ induces an
orthogonal product
\begin{equation}\label{orthprod}
\times:\Lambda^1\otimes \Delta_+\cong \Lambda^1\otimes\Lambda^1\to
\Delta_-\cong\Lambda^1,
\end{equation}
a fact previously noticed in~\cite{hi01}.

Next we analyse the case of $SU(2)\cdot SU(2)\times U(1)$. As before, we label irreducible representations by their highest weight expressed in the basis of fundamental roots. Recall that the irreducible representations of $SU(2)$ are given by the symmetric power $\sigma^n=\odot^n\C^2$ of the complex vector representation $\C^2$ and are labeled by the half--integer $l=n/2$. They are real for $n$ even and quaternionic for $n$ odd. Consequently, the irreducible representations of
$SU(2)\cdot SU(2)\times U(1)$ can be labeled by
$(l_1,l_2,m)=(l_1)\otimes(l_2)\otimes(m)$, where the third factor denotes the irreducible $S^1$--representation $S_m:\theta(z)\mapsto
e^{im\theta}\cdot z$ which is one--dimensional and complex. We will use, as we already did in Theorem~\ref{bracket}, the notation from~\cite{sa89} and denote a real module $V$ by $[n_1,\dots,n_l]$ if its complexification $V\otimes\C=(n_1,\ldots,n_l)$ is self--dual (that is, $V\otimes\C$ is complex irreducible). Otherwise, we write $\llbracket n_1,\ldots,n_l\rrbracket$, which means that $V\otimes\C=W\oplus\overline{W}$ with $W$ an irreducible complex module non--equivalent to $\overline{W}$. By assumption, we have
$$
\Lambda^1=\mf{su}(2)\oplus\mf{su}(2)\oplus
\R^2=[1,0,0]\oplus[0,1,0]\oplus\llbracket0,0,2\rrbracket.
$$
Hence, $SU(2)\cdot SU(2)\times U(1)$ acts with weights
$0,\,\alpha_1,\,\alpha_2$ and $2m$, with $\alpha_1$ and $\alpha_2$ being the fundamental roots of $SU(2)\times SU(2)$. The $Spin(8)$--weights on $\Delta_{\pm}$ are $(x_1\pm\ldots\pm x_4)/2$ with an even (respectively odd) number of minus signs, where the $x_j$ are the parameters of the standard Cartan sub--algebra of $Spin(8)$. Substituting
$$
x_1=\alpha_1,\,x_2=\alpha_2,\,x_3=2m,\,x_4=0
$$
shows that as an $SU(2)\cdot SU(2)\times U(1)$--space,
$$
\Delta_+=\Delta_-=\llbracket\frac{1}{2},\frac{1}{2},1\rrbracket=\llbracket\C^2\otimes\C^2\otimes
S_1\rrbracket.
$$
In particular, the action of $SU(2)\cdot SU(2)\times U(1)$ on $\Delta_{\pm}$ preserves a
complex structure. Note however that this structure does not reduce to $SU(4)$ as the torus component acts non--trivially on
$\lambda^{4,0}\Delta_{\pm}$. Permuting with the triality automorphisms yields a complex structure on $\Lambda^1$ and $\Delta_{\pm}$ if the isometry is an element of $\Lambda^3\Delta_{\mp}$.
 
Finally we consider the group $Sp(1)\cdot Sp(2)$, that is
$$
\Lambda^1=\mf{su}(2)\oplus \mf{z}^5=[1,0,0]\oplus[0,2,-1].
$$
Here the first component refers to the representation labeled by
$\alpha$, the fundamental root of $\mf{sp}(1)\otimes\C=\mf{su}(2)\otimes\C$, while the last two indices $(m_1,m_2)$ designate the irreducible
$Sp(2)$--representation with respect to the basis of fundamental
roots $\beta_1$ and $\beta_2$. The weights of the action on
$\Lambda^1$ are
$0,\,\alpha,\,\beta_1+\frac{1}{2}\beta_2,\,\beta_1+\frac{3}{2}\beta_2$.
Substituting as above, we obtain
$$
\Delta_+=\Delta_-=[1/2,1,1]=[\C^2\otimes\H^2],
$$
where the quaternionic space $\H^2$ serves as a model for the irreducible spin
representation of $Sp(2)=Spin(5)$.

In this paper we will focus on geometric structures associated with the groups $PSU(3)$ and $Sp(1)\cdot Sp(2)$ stabilising a supersymmetric map $\sigma_{\pm}\in\Lambda^3\Delta_{\pm}\subset\Delta_-\otimes\Lambda^1$, thus acting irreducibly on $\Lambda^1$. Before we continue, a thorough discussion of the linear algebra of these groups is in order. 

We begin with the group $PSU(3)=SU(3)/\ker Ad$ whose (negative definite) Killing form we denote by $B$. It is the identity component of the automorphism group of $\mf{su}(3)$ and therefore compact and of dimension $8$. In particular, the adjoint representation $Ad:SU(3)\to
SO(8)$ descends to an embedding $PSU(3)\hookrightarrow SO(8)\subset GL(8)$, so that
$\Lambda^1=\mf{su}(3)$. The group $PSU(3)$ arises as the stabiliser of the 3--form
$$
\rho(X,Y,Z)=-\frac{4}{3}B([X,Y],Z)
$$
inside $GL_+(8)$, the linear transformations of positive determinant. Further, the $GL_+(8)$--orbit of $PSU(3)$--invariant forms is {\em open}, i.e. they are {\em stable} following the language of~\cite{hi01}, for $\dim GL_+(8)-\dim PSU(3)=\dim\Lambda^3$. As we have already used above, a $PSU(3)$--invariant form $\rho$ can be expressed in a {\em $PSU(3)$--frame} as
\begin{equation}\label{psu3form}
\rho=\frac{1}{2}e_{123}+\frac{1}{4}e_1(e_{47}-e_{56})+\frac{1}{4}e_2(e_{46}+e_{57})+\frac{1}{4}e_3(e_{45}-e_{67})+\frac{\sqrt{3}}{4}e_8(e_{45}+e_{67}).
\end{equation}
When dealing with $PSU(3)$, we always assume to work with such a frame unless otherwise stated.

Next we will discuss some elements of the representation theory for $PSU(3)$. The Lie algebra of the stabiliser of $\rho$ inside $\Lambda^2$ is given by the vectors $x\llcorner\rho$, $x\in\Lambda^1$. A maximal torus is spanned by $x_3=e_3\llcorner\rho$ and $x_8=e_8\llcorner\rho$ with roots $\pm\alpha_1=\pm i(x^3+\sqrt{3}x^8)/4$, $\pm\alpha_2=\pm i(x^3-\sqrt{3}x^8)/4$ and $\pm(\alpha_1+\alpha_2)=\pm ix^3/2$ and root vectors $x_{\alpha_1} = e_4-ie_5$, $x_{\alpha_2} = e_6+ie_7$ and
$x_{\alpha_1+\alpha_2} = e_1-ie_2$. For the exterior algebra we find the following decomposition, where $\Lambda^p_q$ represents a $q$--dimensional irreducible subspace of $\Lambda^p$.

\begin{prop}\label{psu3decomp}\hfill\newline\vspace{-0.8cm}
\begin{enumerate}
\item $\Lambda^1=\mf{su}(3)=[1,1]$ is irreducible.
\item
$\Lambda^2=\Lambda^2_8\oplus\Lambda^2_{20}=[1,1]\oplus\llbracket1,2\rrbracket$. 
\item
$\Lambda^3=\Lambda^3_1\oplus\Lambda^3_8\oplus\Lambda^3_{20}\oplus\Lambda^3_{27}=\mb{1}\oplus[1,1]\oplus\llbracket1,2\rrbracket\oplus[2,2]$
\item
$\Lambda^4=2\Lambda^4_8\oplus2\Lambda^4_{27}=2[1,1]\oplus2[2,2]$
\end{enumerate}
\end{prop}

Note that the Hodge $\star$--operator equivariantly identifies $\Lambda^p$ with $\Lambda^{8-p}$. This decomposition can be also understood from a cohomological point of view well--suited for later purposes. The $\mf{su}(3)$--structure on $\Lambda^1$ induces a $PSU(3)$--invariant operator $c_k:\Lambda^k\to\Lambda^{k+1}$ by extension of
$$
c_1e_i=\sum_{j<k}c_{ijk}e_j\wedge e_k
$$
built out of the structure constants $c_{ijk}$ of $\mf{su}(3)$. It is therefore just the exterior differential operator restricted to the left--invariant differential forms of $SU(3)$ with adjoint $c^*=d^*=-\star\,d\,\star$. The resulting elliptic complex is isomorphic to the de Rham cohomology $H^*\big(SU(3),\R\big)$ which is trivial except for the Betti numbers $b^0=b^3=1=b^5=b^8$. Hence, $\im c_k=\ker c_{k+1}$ for $k=0,\,1,\,3,\,5,\,6$ and $\im c_k=\ker c_{k+1}\oplus \R$ for $k=-1,\,2,\,4,\,7$. Schematically, we have
\begin{equation}\label{b_decomp}
\begin{array}{lclclclclclclclcl}
\Lambda^0_1 & \phantom{\stackrel{b}{\to}} & & & & & \Lambda^3_1 & & & & \Lambda^5_1 & & & & & \phantom{\stackrel{c}{\to}} & \Lambda^8_1\\[2pt]
& & \Lambda^1_8 & \stackrel{c}{\to} & \Lambda^2_8 & & \Lambda^3_8 & \stackrel{c}{\to} & \Lambda^4_8 & & & & \Lambda^6_8 & \stackrel{c}{\to} & \Lambda^7_8 & & \\[2pt]
& & & & & &  & & \Lambda^4_8 &  \stackrel{c}{\to} & \Lambda^5_8 & & & & & & \\ 
& & & & \Lambda^2_{20} & \stackrel{c}{\to} & \Lambda^3_{20} & & & & \Lambda^5_{20} & \stackrel{c}{\to} & \Lambda^6_{20} & & & & \\[2pt]
& & & & & & \Lambda^3_{27} & \stackrel{c}{\to} & \Lambda^4_{27} & & & & & & & & \\[2pt]
& & & & & & & & \Lambda^4_{27} & \stackrel{c}{\to} & \Lambda^5_{27} & & & & &
\end{array}
\end{equation}
with an arrow indicating the non--trivial maps. In particular, we will use the more natural splitting of $\Lambda^4$ into $\Lambda^4_o=\ker c_3$ and $\Lambda^4_i=\im c^*_5$ instead of the $SO(8)$--equivariant splitting into self-- and anti--self--dual forms. From this, we can also construct the projection operators for $\Lambda^2\otimes\C$, which will be useful in Section~\ref{Dirac}.

\begin{prop}\label{liealgsu3} For any $\alpha\in\Lambda^2$ we have $c_2(\alpha)=-\alpha^*\rho$. Moreover, $\Lambda^2_8=\ker c_2$ and the projection operator on the complement is $\pi^2_{20}(\alpha)=\frac{4}{3}c^*_3c_2(\alpha)$. For the complexification, we find $\Lambda^2_{20}\otimes\C=\Lambda^2_{10+}\oplus\Lambda^2_{10-}=(1,2)\oplus (2,1)$,
where
$$
\Lambda^2_{10\pm}=\{\alpha\in\Lambda^2\otimes\C\,|\,\star(\rho\wedge \alpha)=\pm i\sqrt{3}\alpha^*\rho\}.
$$
The projection operators are $\pi^2_{10\pm}(\alpha)= \frac{2}{3}c^*_3c_2(\alpha)\mp\frac{8\sqrt{3}}{9}i\star(c_2(\alpha)\wedge\rho)$.
\end{prop}

The proof can be readily verified by applying Schur's Lemma with the sample vectors $x_{\alpha_2}\wedge x_{\alpha_1+\alpha_2}\in(1,2)$ and $x_{\alpha_1}\wedge x_{\alpha_1+\alpha_2}\in(2,1)$.

The $PSU(3)$--invariant supersymmetric maps $\sigma_{\pm}$ in $\ker\mu_{\pm}\cong\Lambda^3\Delta_{\mp}$ are characterised (up to a scalar) by the equations
$$
x\llcorner\rho(\sigma_{\pm})\,=\,\frac{1}{2}\kappa(x\llcorner\rho)\cdot\sigma_{\pm}-\sigma_{\pm}\circ x\llcorner\rho\,=\,0.
$$
Their matrices with respect to a $PSU(3)$--frame and a fixed orthonormal basis of $\Delta_{\pm}$ are given in Appendix~\ref{susymatrix}. Note that their determinant is $1$ since the outer triality morphisms reverse the orientation.

We close our discussion of $PSU(3)$ with a remark on special $PSU(3)$--orbits in the Grassmannians $\widetilde{G}_3(\Lambda^1)$ and $\widetilde{G}_5(\Lambda^1)$ of oriented 3-- and 5--dimensional planes in $\Lambda^1$. These orbits consist of calibrated planes, a notion due to Harvey and Lawson~\cite{hala82} which we briefly recall. Let $(V,g,\tau)$ be an oriented (real) vector space with a Euclidean metric $g$ and a $k$--form $\tau\in\Lambda^kV^*$. We say that $\tau$ defines a {\em
calibration} if for every oriented $k$-plane $\xi=f_1\wedge\ldots\wedge f_k$ in $V$ given by some orthonormal system $f_1,\ldots,f_k$, the inequality $\tau(f_1,\ldots,f_k)\le 1$
holds and is met for at least one $k$--plane. Such a plane is said to be {\em calibrated} by $\tau$. A classical example is provided by the imaginary octonions whose so--called associative and co--associative planes are calibrated by the $G_2$--invariant forms $\varphi$ and $\star\varphi$ respectively.

\begin{prop}
Let $\rho$ be the $PSU(3)$--invariant $3$--form~(\ref{psu3form}) and $\tau=2\rho$. Then $\tau(\xi)\le 1$ with equality if and only if
$\xi=Ad(A)\mf{h}$ for $A\in SU(3)$, where $\mf{h}$ is a suitably oriented $\mf{su}(2)$--subalgebra associated with a highest root. Furthermore,
$\star\tau(\xi)\le 1$ with equality if and only if $\xi$ is perpendicular to a 3--plane calibrated by $\tau$. In particular,
$PSU(3)$ acts transitively on the set of calibrated $3$-- and $5$--planes.
\end{prop}

\begin{pf}
We adapt the proof from~\cite{ta85}. Let $e_1,\ldots,e_8$ be a
$PSU(3)$--frame inducing the Euclidean norm $\parallel\cdot\parallel$, and fix
the Cartan subalgebra $\mf{t}$ spanned by $e_3$ and $e_8$. Let $E_1=e_5, F_1=-e_4$, $E_2=-e_6, F_2=e_7$ and $E_3=e_1, F_3=e_2$, and put $\lambda_1=(e^3+\sqrt{3}e^8)/4$, $\lambda_2=(e^3-\sqrt{3}e^8)/4$ and $\lambda_3=\lambda_1+\lambda_2=e^3/2$. Then $\parallel\lambda_i\parallel=1/2$ and we immediately verify the relations
\begin{equation}\label{su2rel}
[T,E_i]=\lambda_i(T)F_i,\quad
[T,F_i]=-\lambda_i(T)E_i\quad\mbox{and}\quad
[E_i,F_i]=\lambda_i
\end{equation}
for $T\in\mf{t}$ and $i=1,2,3$. Next let $\xi\in\widetilde{G}_3(\Lambda^1)$.
Since $\mf{t}$ is a Cartan subalgebra, $Ad\big(SU(3)\big)X\cap\mf{t}\not=\emptyset$ for any $0\not=X\in\Lambda^1$.
Moreover, $\rho$ is $Ad$--invariant, so we may assume that $\xi\cap\mf{t}\not=\emptyset$ up to the action of an element in $SU(3)$.
Pick $T\in\xi\cap\mf{t}$ and extend it to a positively oriented basis $\{T,X,Y\}$ of $\xi$. Then
$$
X=T_0+\sum\limits_{i=1}^3s_iE_i+\sum\limits_{i=3}^3t_iF_i,
$$
where $T_0\in\mf{t}$. Hence $\parallel X\parallel^2=\parallel T_0\parallel^2+\sum_{i=1}^3(s^2_i+t^2_i)$ which implies
\begin{eqnarray}
\parallel[T,X]\parallel^2 & = & \parallel\sum\limits_{i=1}^3s_i\lambda_i(T)F_i-t_i\lambda_i(T)E_i\parallel^2\nonumber\\
& \le & \sum\limits_{i=1}^3(s_i^2+t^2_i)\parallel\lambda_i\parallel^2\parallel T\parallel^2\nonumber\\
& \le & \frac{1}{4}\parallel T\parallel^2\parallel X\parallel^2\label{calineq}.
\end{eqnarray}
The Cauchy--Schwarz inequality yields
\begin{equation}\label{csineq}
|\tau(T,X,Y)|=2|g([T,X],Y)|\le2\parallel[T,X]\parallel\parallel Y\parallel\le\parallel T\parallel\parallel X\parallel\parallel Y\parallel,
\end{equation}
so that $\tau(f_1,f_2,f_3)\leq1$ for any plane $\xi=f_1\wedge f_2\wedge f_3$. Furthermore, equality holds for~(\ref{csineq}) if and only if $Y$ is a multiple of $[T,X]$ and $\tau(T,X,Y)>0$. For~(\ref{calineq}), equality holds if and only if $T\in\R\lambda_i$ and
$X\in\langle E_i, F_i \rangle$ for an $i\in\{1,2,3\}$. Consequently, if $\xi=\langle T,X,Y\rangle$ is calibrated, then
$Y\in\langle E_i, F_i \rangle$ and because of
(\ref{su2rel}), $\xi$ is an $\mf{su}(2)$--algebra.

Since $(\star\tau)_{|\xi^{\perp}}=\star(\tau_{|\xi})$ any
calibrated 5--plane is the orthogonal complement of an $\mf{su}(2)$--algebra. Moreover, any two subalgebras of highest root are conjugate.
\end{pf}

\begin{rem}
As for $G_2$-- or $Spin(7)$--structures, calibrations give rise to a natural type of submanifolds for $PSU(3)$--structures, namely those whose tangent space at any point is calibrated. More generally, Tasaki showed~{\rm \cite{ta85}} that for any compact simple Lie group $G$ with Killing form $B$ and Lie algebra $\mathfrak{g}$, the 3--form
$$
\tau(X,Y,Z)=-\frac{1}{\parallel\delta\parallel}B([X,Y],Z),
$$
where $\parallel\delta\parallel$ is the norm with respect to $B$ of a highest root $\delta$ of $\mathfrak{g}$, defines a calibration on $G$. Furthermore, any calibrated submanifold is a translate of a compact simple $3$--dimensional subgroup associated with $\delta$.
\end{rem}

Next we turn to the group $Sp(1)\cdot Sp(2)$. Here the vector representation of $GL(8)$ restricted to this group is $\Lambda^1=[\C^2\otimes\H^2]$. Elevating this to the fourth exterior power yields an invariant 4--form $\Omega$. To describe $\Omega$ explicitly, think of $\Lambda^1$ as a quaternionic vector space $\O\cong\H^2$. This is acted on by $Sp(2)$ which fixes the three K\"ahler 2--forms $\omega_i$, $\omega_j$ and $\omega_k$ given by $\omega_i(x,y)=g(x\cdot i,y)$ etc. Then $\Omega=\omega_i\wedge\omega_i+\omega_j\wedge\omega_j+\omega_k\wedge\omega_k$ is $Sp(1)\cdot Sp(2)$--invariant~\cite{kr66}. In terms of the orthonormal basis~(\ref{octbasis}), we find $\omega_i = e_{12}-e_{34}+e_{56}-e_{78}$, $\omega_j = e_{13}+e_{24}+e_{57}+e_{68}$ and
$\omega_k = e_{14}-e_{23}+e_{58}-e_{67}$,
so that
\begin{eqnarray}
\frac{1}{2}\Omega & = & -3e_{1234}+e_{1256}-e_{1278}+e_{1357}+e_{1368}+e_{1458}-e_{1467}\nonumber\\
& & -e_{2358}+e_{2367}+e_{2457}+e_{2468}-e_{3456}+e_{3478}-3e_{5678}\label{sp1sp2form}.
\end{eqnarray}
In analogy with the $PSU(3)$--case we refer to any orthonormal frame $e_1,\ldots,e_8$ such that $\Omega/2$ is of the form~(\ref{sp1sp2form}) as an $Sp(1)\cdot Sp(2)$--{\em frame}.

The invariant 4--form induces a splitting of $\mf{so}(8)$ into the Lie algebra of the stabiliser and its orthogonal complement which we need to make explicit. If $a^*\Omega=0$ for $\sum_{i<j}a_{ij}e_i\wedge e_j$, where $a^*$ denotes the usual action of $\mf{gl}(8)$ on exterior forms, then
$$
\begin{array}{lll}
a_{68} - a_{13} - a_{24} + a_{57}=0,\,\, & a_{46} - a_{17}=0,\,\, & a_{47} - a_{25}=0, \\
a_{23} - a_{14}-a_{67}+a_{58}=0,\,\, & a_{35} +a_{17}=0,\,\, & a_{28} + a_{17}=0,
\\
a_{34} - a_{78}+a_{56}-a_{12}=0,\,\, & a_{45} + a_{18}=0,\,\, & a_{26} - a_{48} = 0,\\
a_{38} + a_{25}=0,\,\, & a_{16} + a_{25}=0\,\, & a_{27} - a_{18}=0,\\
a_{36} + a_{18}=0,\,\, & a_{15} - a_{48}=0,\,\, & a_{37} - a_{48}=0
.
\end{array} 
$$
A maximal torus is spanned for instance by $a_1=(e_{12}-e_{34}+e_{56}-e_{78})/2$, $a_2=e_{12}+e_{34}+e_{56}+e_{78}$ and $a_3=e_{12}+e_{34}-e_{56}-e_{78}$ with 
corresponding fundamental roots $\alpha=ia^1$, $\beta_1=2i(a^2-a^3)$ and $\beta_2=2ia^5$. The weights of $\Lambda^1=[\C^2\otimes\H^2]$ are
\begin{equation}\label{weightssp1sp2}
\pm{\textstyle\frac{1}{2}}(\alpha+\beta_1),\quad 
\pm{\textstyle\frac{1}{2}}(\alpha-\beta_1),\quad 
\pm{\textstyle\frac{1}{2}}(\alpha+\beta_1)+\beta_2,\quad
\pm{\textstyle\frac{1}{2}}(\alpha-\beta_1)-\beta_2
\end{equation}
and the weight vectors are given by $x_{(\alpha+\beta_1)/2}=e_5-ie_6$, $x_{(\alpha-\beta_1)/2}=e_7+ie_8$,  $x_{(\alpha+\beta_1)/2+\beta_2}=e_1-ie_2$ and $x_{(\alpha-\beta_1)/2-\beta_2}=e_3+ie_4$. 
A different characterisation of the decomposition $\Lambda^2=\mf{sp}(1)\oplus\mf{sp}(2)\oplus\big(\mf{sp}(1)\oplus\mf{sp}(2)\big)^{\perp}$ is given by the equivariant map $\alpha\mapsto\alpha\llcorner\Omega$. A straightforward application of Schur's Lemma yields

\begin{prop}\label{2formssp1sp2} We have
$\mf{sp}(1)=\{\alpha\in\Lambda^2\,|\,(\alpha\llcorner\Omega)=5\alpha\}$, $\mf{sp}(2)=\{\alpha\in\Lambda^2\,|\,(\alpha\llcorner\Omega)=-3\alpha\}$ and $\big(\mf{sp}(1)\oplus\mf{sp}(2)\big)^{\perp}=\{\alpha\in\Lambda^2\,|\,(\alpha\llcorner\Omega)=\alpha\}$. Further, the projection operators onto these modules are
$\pi^2_3(\alpha)=\big(-3\alpha+2\alpha\llcorner\Omega+(\alpha\llcorner\Omega)\llcorner\Omega\big)/32$, $\pi^2_{10}(\alpha)=\big(5\alpha-6\alpha\llcorner\Omega+(\alpha\llcorner\Omega)\llcorner\Omega\big)/32$ and 
$\pi^2_{15}(\alpha)=\big(15\alpha+2\alpha\llcorner\Omega-(\alpha\llcorner\Omega)\llcorner\Omega\big)/16$ respectively.
\end{prop}

The decomposition on the remaining exterior powers is this. If $\sigma=(1/2)=\C^2$ denotes as above the vector representation of $Sp(1)=SU(2)$ and $\H^2=[1/2,1]$ the vector representation of $Sp(2)$, we have a Clebsch--Gordan like decomposition~\cite{sw89}
$$
\Lambda^r\cong\Lambda^r[\sigma\otimes\H^2]\cong\bigoplus\limits_{s=0}^{[r/2]}[\sigma^{r-2s}\otimes V^r_s],\quad 0\leq r\leq 8,
$$
where the irreducible $GL(2,\H)$--module $V^r_s$ is the direct sum of the irreducible $Sp(2)$--modules
$$
\lambda^r_s=(\frac{r}{2},\big[\frac{3+s}{2}\big]),\quad 0\leq r-2s\leq k
$$
(note that our choice of a basis differs from~\cite{sa89}). In particular, $\lambda^1_0=(\frac{1}{2},1)=\H^2$. Hence:

\begin{prop}\label{sp(1)sp(2)decomp}\hfill\newline\vspace{-0.8cm}
\begin{enumerate}
\item $\Lambda^1=[\sigma\otimes\lambda^1_0]=[\frac{1}{2},\frac{1}{2},1]$ is irreducible.
\item
$\Lambda^2=[\sigma^2]\oplus[\lambda^2_1]\oplus[\sigma^2\otimes\lambda^2_1]=[1,0,0]\oplus[0,1,2]\oplus[1,1,1]$
\item
$\Lambda^3=[\sigma\otimes\lambda^1_0]\oplus[\sigma\otimes\lambda^3_1]\oplus[\sigma^3\otimes\lambda^1_0]=[\frac{1}{2},\frac{1}{2},1]\oplus[\frac{1}{2},\frac{3}{2},2]\oplus[\frac{3}{2},\frac{1}{2},1]$
\item
$\Lambda^4=\R\oplus[\lambda^2_0]\oplus[\lambda^4_2]\oplus[\sigma^2\otimes\lambda^2_0]\oplus[\sigma^2\otimes\lambda^2_1]\oplus[\sigma^4]=[0,0,0]\oplus[0,1,1]\oplus[0,2,2]\oplus[1,1,1]\oplus[1,1,2]\oplus[2,0,0]$
\end{enumerate}
\end{prop}

The $Sp(1)\cdot Sp(2)$--invariant supersymmetric map $\sigma_+$ in $\ker\mu_+$ is determined (up to a scalar) by the equation
$$
a(\sigma_+)\,=\,\frac{1}{2}\kappa(a)\cdot\sigma_+-\sigma_+\circ a\,=\,0,\quad a\in\mf{sp}(1)\oplus\mf{sp}(2).
$$
An explicit matrix representation is given in Appendix~\ref{susymatrix}. Its determinant is $-1$, in accordance with Theorem~\ref{orbitdecomposition} and the $PSU(3)$--case.

Finally, we obtain again a calibration form by taking a suitable multiple of $\Omega$. In~\cite{ta85}, Tasaki proved the

\begin{prop}
Let $\Omega$ be the $Sp(1)\cdot Sp(2)$--invariant $4$--form~(\ref{sp1sp2form}) and $\tau=\Omega/6$. Then $\tau(\xi)\le 1$ with equality if and only if $\xi$ is a suitably oriented $Sp(1)$--invariant $4$--plane. In particular, $Sp(1)\cdot Sp(2)$ acts transitively on the set of calibrated $4$--planes.
\end{prop}
%
%
%
%
%
\section{Topological reductions to $PSU(3)$}
\label{topred}
\begin{defn}
Let $M^8$ be an $8$--dimensional, smooth manifold. A {\em topological $PSU(3)$-- or $Sp(1)\cdot Sp(2)$--structure} is a reduction from the frame bundle on $M$ to a principal $PSU(3)$-- or $Sp(1)\cdot Sp(2)$--fibre bundle.
\end{defn}

A topological $PSU(3)$--structure is equivalent to the choice of an orientation and the existence of a $3$--form $\rho$ with $\rho_x\in\Lambda^3T_x^*M$ lying in the orbit diffeomorphic to $GL_+(8)/PSU(3)$ for any $x\in M$. Similarly, a topological $Sp(1)\cdot Sp(2)$--structure is tantamount to endowing $M$ with a $4$--form $\Omega$ such that $\Omega_x\in\Lambda^4T_x^*M$ lies in the orbit diffeomorphic to $GL(8)/Sp(1)\cdot Sp(2)$ for all $x\in M$. In this section, we investigate necessary and sufficient criteria for a $PSU(3)$--reduction to exist. For the $Sp(1)\cdot Sp(2)$--case, one has the following result.

\begin{thm}\cf{cava98d} 
Let $M$ be an oriented closed connected spinnable manifold of dimension 8. If $M$ carries an $Sp(1)\cdot Sp(2)$--structure, then $8e+p_1^2-4p_2=0$. Moreover, provided that $H^2(M,\Z_2)=0$, we have $w_6=0$ and there exists an $R\in H^4(M,\Z)$ such that $Sq^2p_2R=0$, $(Rp_1-2R^2)[M]\equiv0\mod 16$ and $(R^2+Rp_1-e)[M]\equiv 0\mod 4$, where $[M]$ denotes the fundamental class of $M$. Conversely, these conditions are sufficient (regardless of $H^2(M,\Z_2)=0$) to ensure the existence of an $Sp(1)\cdot Sp(2)$--structure.
\end{thm}

Necessary conditions for a topological $PSU(3)$--structure to exist easily follow from a characteristic class computation using the Borel--Hirzebruch formalism~\cite{bohi58}. Let $\pm x_1,\ldots,\pm x_4$ denote the weights of the vector representation of $SO(8)$. Formally, the total Pontrjagin class $p$ and the
Euler class $e$ of $M$ are expressed as the product
$$
p=\prod(1+x_i^2),\quad e=\prod x_i.
$$
If the tangent space is associated with the adjoint representation of $PSU(3)$, the $SO(8)$--weights become the $PSU(3)$--roots under restriction. Substituting 
$$
x_1=\alpha, \quad x_2=\beta, \quad x_3=\alpha+\beta, \quad x_4=0,
$$ 
a reduction to $PSU(3)$ implies $p_1=2(\alpha^2+\alpha\beta+\beta^2)$ and
$p_2=\alpha^4+2\alpha^3\beta+3\alpha^2\beta^2+2\beta^3\alpha+\beta^4$,
hence $4p_2=p_1^2$. Moreover, we obviously have $e=0$. A first consequence is the following classification result.

\begin{prop}\label{topclass}
Let $(G/H,g)$ be a compact Riemannian homogeneous space with $G$ simple. If $M=G/H$ admits a topological reduction to $PSU(3)$, then $G/H$ is diffeomorphic to $SU(3)$.
\end{prop}

\begin{pf}
Since $G$ sits inside the isometry group of $(M,g)$, its dimension
is less than or equal to $9\cdot 8/2=36$. If we had equality, then
$M$ would be diffeomorphic to a torus or, up to a finite covering, to an $8$--sphere. While the first case is ruled out for $G$ has to
be simple, the second case is excluded since $e(S^8)\not= 0$.
Hence $G$ must be, up to a covering, a group of type
$A_1,\ldots,A_5$, $B_2,B_3$, $C_3$, $D_4$ or $G_2$. As a closed
subgroup of $G$, $H$ is compact and hence reductive. Therefore $H$
is covered by a direct product of simple Lie groups and a torus,
that is the Lie algebra of $H$ is isomorphic to
$\mf{h}=\mf{g}_1\oplus\ldots\oplus\mf{g}_k\oplus\mf{t}^l$. If we
denote by $rk(G)$ the rank of the Lie group $G$, we get the following
necessary conditions:
$$ \begin{array}{l}
k\le rk(G) \\
l+\sum rk(\mf{g}_i)\le rk(G)\\
l+\sum{dim(\mf{g}_i)}=dim(G)-8,
\end{array}
$$
which yields the possibilities displayed in the table below.
\begin{table}[hbt]
\begin{center}
\begin{tabular}{c|c|c|c}
$G$ & $H$ up to a covering & $dim(H)$ & $rk(H)$\\
\hline
$A_2$ & $\{1\}$ & $0$ & $0$\\
$A_3$ & $A_1\times A_1\times S^1$ & $7$ & $3$\\
$A_4$ & $A_3\times S^1,\,\, G_2\times S^1\times S^1,\,\, A_2\times A_2$ & 16 & 4\\
$A_5$ & $A_1\times A_4,\,\, A_1\times A_1\times B_3,\,\, A_1\times
A_1\times C_3$ & $27$ & $5$\\
$B_2$ & $S^1\times S^1$ & $2$ & $2$\\
$B_3$ & $A_1\times B_2$ & $13$ & $3$\\
$C_3$ & $A_1\times B_2$ & $13$ & $3$\\
$D_4$ & $A_1\times A_1\times G_2$ & $20$ & $4$\\
$G_2$ & $A_1\times A_1$ & $6$ & $2$\\
\end{tabular}
\end{center}
\end{table}
\newline It follows that $H$ is of maximal rank, that is $rk(H)=rk(G)$,
unless $G=SU(3)$ and $H=\{1\}$. But in the first case,~\cite{sa58} implies $e(G/H)\not=0$, whence the assertion.
\end{pf}

Since $\pi_1\big(PSU(3)\big)=\Z_3$, the inclusion $PSU(3)\subset SO(8)$ lifts to $Spin(8)$. In particular, any $8$--manifold admitting a topological $PSU(3)$--structure must be spinnable, hence the first and second Stiefel--Whitney class $w_1$ and $w_2$ of $M$ have to vanish. By a straightforward computation using the definition of the $\widehat{A}$--genus and the signature of $M$, $sgn(M)=b_4^+-b_4^-$, where $(b_4^+,b_4^-)$ is the signature of the Poincar\'e pairing on $H^4(M,\Z)$, we derive the following

\begin{lem}
If $M$ is a compact spin manifold with $p_1^2=4p_2$, then
$sgn(M)=16\widehat{A}[M]$. In particular, $sgn(M)\equiv 0\,\,\mod 16$.
\end{lem}

\begin{cor}
Let $M$ be a compact simply--connected manifold whose frame bundle reduces to $PSU(3)$. If $\widehat{A}[M]=0$ (e.g. if there exists a metric
with strictly positive scalar curvature), then $1+b_2+b_4^+=b_3$.
\end{cor}

\begin{exmp}
As already stated in Section~\ref{triality}, the Betti numbers $b^q$ of $SU(3)$ are either $0$ or $1$ for $q=0$, $3$, $5$ or $8$, in accordance with the corollary.
\end{exmp}

As $e=0$ and $sgn(M)\equiv 0\mod4$, we can assert the existence of two linearly independent vector fields~\cite{th69}. The orthogonal product $\times$ in (\ref{orthprod}) produces a third one. In particular, $w_6=0$. Taking $k=0$ in the following proposition yields the existence of four pointwise linearly independent vector fields.

\begin{prop}\cf{cava98a}
Let $M$ be a closed connected smooth spin manifold of dimension 8. If $w_6=0$, $e=0$ and $\{4p_2-p_1^2\}[M]\equiv 0 \mod 128$, and if there is a $k\in\Z$ such that $4p_2=(2k-1)^2p_1^2$ and $k(k+2)p_2[M]\equiv 0 \mod 3$, then $M$ has four pointwise linearly independent vector fields.
\end{prop}

As a further consequence, $w_5=0$.

\begin{prop}\label{wu}
We have $w_4^2=0$. In particular, all Stiefel-Whitney numbers vanish.
\end{prop}

\begin{pf}
By Wu's formula,
$$
{\rm
Sq}^k(w_m)=w_kw_m+\binom{k-m}{1}w_{k-1}w_{m+1}+\ldots+\binom{k-m}{k}w_0w_{m+k}.
$$
A further theorem of Wu asserts that
$$
w_k=\sum_{i+j=k}{\rm Sq}^i(v_j),
$$
where the elements $v_k\in H^k(M,\Z_2)$ are defined through the identity $v_k\cup x[M]={\rm Sq}^k(x)[M]$ for $x\in
H^{n-k}(M,\Z_2)$. In particular, we have $v_i=0$ for $i>4$. It follows that $v_1=v_2=v_3=0$, $w_4=v_4$ and $w_8={\rm
Sq}^4w_4=w^2_4=0$.
\end{pf}

We summarise our results in the following proposition.

\begin{prop}\label{obpsu3}
If a closed and oriented $8$--manifold $M$ carries a topological $PSU(3)$--structure, then all Stiefel--Whitney classes vanish except $w_4$, and $w_4^2=0$. Moreover, we have $e=0$ and $p_1^2=4p_2$. There exist four linearly independent vector fields on $M$ and all Stiefel-Whitney numbers vanish.
\end{prop}

Finding sufficient conditions to ensure the existence of a topological $PSU(3)$--structure over closed $M$ occupies us next. This problem is considerably harder than the analogous problem for topological $G_2$--structures on $7$--manifolds. Here, a reduction to $G_2$ implies that the underlying manifold is spin. Conversely, assuming that it is spin, we can pick a spin structure and consider the associated spinor bundle $\Delta$. This is a real bundle of rank $8$ whose sphere bundle is associated with $Spin(7)/G_2$. The existence of a topological $G_2$--structure is therefore equivalent to the existence of a nowhere vanishing spinor field, for which there is no obstruction since the Euler class of $\Delta$ vanishes trivially on dimensional grounds. Similarly, we deduce from Proposition~\ref{obpsu3} and the discussion in Section~\ref{triality} that a topological $PSU(3)$--structure on a Riemannian $8$--manifold can be characterised by a spinor--valued $1$--form. However, this must be a section of special algebraic type and taking an arbitrary, nowhere vanishing section will in general not result in a topological reduction to $PSU(3)$. 

Therefore we restrict ourselves to a special class of $PSU(3)$--structures for which the problem of finding sufficient conditions becomes easier. Assume that $M$ admits a principal $SU(3)$--fibre bundle $\widetilde{P}$ such that the adjoint bundle $\mf{su}(\widetilde{P})=\widetilde{P}\times_{SU(3)}\mf{su}(3)$
is isomorphic with the tangent bundle. Then $TM$ is naturally associated with the $PSU(3)$--structure $P=\widetilde{P}/\Z_3$ induced by the exact sequence (where $\Z_3$
is central)
$$
1\to\Z_3\to SU(3)\stackrel{p}{\to} PSU(3)\to 1.
$$
Clearly, not every $PSU(3)$--structure arises this way: The set ${\rm Prin}_G(M)$ of principal $G$--fibre bundles over $M$ can be identified with $H^1(M,G)$ (e.g~\cite{lami89} Appendix A). The sequence above gives rise to the exact sequence
$$
\ldots\to H^1(M,\Z_3)\to {\rm Prin}_{SU(3)}(M) \stackrel{p_*}{\to}
{\rm Prin}_{PSU(3)}(M)\stackrel{t}{\to} H^2(M,\Z_3).
$$
Hence, a principal $PSU(3)$--bundle $P$ is induced by an $SU(3)$--bundle if and only if the obstruction class $t(P)\in H^2(M,\Z_3)$ vanishes. Following~\cite{avis79}, we call this class the {\em triality class}. By the universal coefficients theorem this obstruction vanishes trivially if $H^2(M,\Z)=0$ and $H^3(M,\Z)$ has no torsion elements of order divisible by three. If $f:M\to BPSU(3)$ is a classifying map for $P$, then $t(P)=f^*t$ for the {\em universal triality class} $t\in H^2(BPSU(3),\Z_3)$. It
is induced by $c_1(E_{U(3)})$, the first Chern class of the universal $U(3)$--bundle $E_{U(3)}$~\cite{wo82a}. Concretely, let
$\overline{p}:U(3)\to PU(3)$ denote the natural projection. The inclusion $SU(3)\subset U(3)$ induces an isomorphism between $PSU(3)$ and $PU(3)$ and therefore identifies $BPSU(3)$ with
$BPU(3)$. Since $BPU(3)$ is simply--connected and
$\pi_2\big(BU(3)\big)=\Z\to\pi_2\big(BPU(3)\big)=\Z_3$ is the reduction $\mod3$ map $\rho_3:\Z\to\Z_3$, the Hurewicz isomorphism theorem and the universal coefficients theorem imply that $B\overline{p}^*:H^2(BPU(3),\Z_3)\to H^2(BU(3),\Z_3)$ is an isomorphism and $H^2(BPU(3),\Z_3)=\Z_3$. Then
$$
t=(B\overline{p}^*)^{-1}\rho_{3*}c_1(E_{U(3)}).
$$
Finding conditions ensuring the existence of topological $PSU(3)$--structures with vanishing triality class therefore boils down to finding conditions for principal $SU(3)$--fibre bundles $\widetilde{P}$ with $\mf{su}(\widetilde{P})\cong TM$.

\begin{thm}\label{psu3structure}
Suppose that $M$ is a connected and closed spin manifold of dimension $8$. Then $TM\cong\mf{su}(\widetilde{P})$ for some principal $SU(3)$--bundle $\widetilde{P}$ if and only if $e=0$, $4p_2=p_1^2$, $w_6=0$, $p_1$ is divisible by 6 and $p_1^2[M]\in 216\Z$.
\end{thm}

\begin{pf}
Let us start with the necessity of the conditions. Since $\mathfrak{su}(3)\otimes\C=\mathfrak{sl}(3,\C)$, the complexification $TM\otimes\C$ equals $\End_0(E)$, the bundle of trace--free complex endomorphisms of $E=\widetilde{P}\times_{SU(3)}\C^3$. The Chern character of $TM\otimes\C$ equals
$$
ch(TM\otimes\C)=8+p_1+\frac{1}{12}(p_1^2-2p_2)
$$
(see for instance~\cite{sa82}). On the other hand,
$$
ch\big(End(E)\big)=ch(E\otimes\overline{E})=1+ch\big(End_0(E)\big).
$$
Now for a complex vector bundle with $c_1(E)=0$,
$$
ch(E)=3-c_2(E)+\frac{1}{2}c_3(E)+\frac{1}{12}c_2(E)^2
$$
and $c_i(E)=(-1)^ic_i(\overline{E})$, which implies
$$
ch(E\otimes\overline{E})=ch(E)\cup
ch(\overline{E})=9-6c_2(E)+\frac{3}{2}c_2(E)^2.
$$
As a consequence,
\begin{equation}\label{pontrjagin}
p_1=p_1\big(\mf{su}(E)\big) = -6c_2(E),\quad p_2=p_2\big(\mf{su}(E)\big)= 9c_2(E)^2.
\end{equation}
In particular, $p_1$ is divisible by 6 and we also rederive the relation $4p_2=p_1^2$. Moreover, $M$ is spinnable, hence the spin
index $\widehat{A}\cup{\rm ch}(E)[M]$ is an integer. Since
$$
\widehat{A}=1-\frac{1}{24}p_1+\frac{1}{5760}(-4p_2+7p_1^2)=1-\frac{1}{24}p_1+\frac{1}{960}p_1^2,
$$
it follows
$$
\widehat{A}\cup
ch(E)[M]=3\widehat{A}[M]+\frac{1}{24}p_1c_2(E)+\frac{1}{12}c_2(E)^2=3\widehat{A}[M]-\frac{1}{216}p_1^2[M]\in\Z.
$$
This means $p_1^2[M]\in 216\Z$, proving the necessity of the conditions.

For the converse I am indebted to ideas of M.~Crabb. Let $B\subset M$ be an embedded open disc in
$M$ and consider the exact sequence of K--groups
$$
K(M,M-B)\to K(M)\to K(M-B).
$$
We have $K(M,M-B)=\widetilde{K}(S^8)\cong\Z$ and the sequence is
split by the spin index
$$
x\in K(M)\mapsto\widehat{A}\cup{\rm ch}(x)[M]\in\Z,
$$
which therefore classifies the stable extensions over $M-B$ to $M$.
The first step consists in finding a stable complex vector bundle
$\xi$ over $M-B$ such that the associated adjoint bundle $\mf{su}(\xi)$ is stably
equivalent to $TM_{|M-B}$ and $c_1(\xi)=0$. To that end, let
$[(M-B)_+,BSU(\infty)]\subset K(M-B)$ denote the set of pointed
homotopy classes, the subscript $+$ indicating a disjoint basepoint.
Let $(c_2,c_3)$ be the map which takes an equivalence class of
$[(M-B)_+,BSU(\infty)]$ to the second and third Chern class of the
associated bundle.

\begin{lem}\label{steenrod}
The image of the map
$$
(c_2,c_3):[(M-B)_+,BSU(\infty)]\to H^4(M,\Z)\oplus H^6(M,\Z)
$$
is the set $\{(u,v)\,\,|\,\,Sq^2\rho_2u=\rho_2v\}$, where $\rho_2:\Z\to\Z_2$ is reduction $\mod2$.
\end{lem}

\begin{pf}
We first prove that for a complex vector bundle $\xi$ with $c_1(\xi)=0$, we have
\begin{equation}\label{square}
Sq^2\rho_2c_2(\xi)=\rho_2c_3(\xi).
\end{equation}
Now if $W_i$ denote the Stiefel--Whitney classes of the {\em real} vector bundle underlying $\xi$, this is equivalent to $Sq^2W_4=W_6$. On the other hand, Wu's formula implies
$$
Sq^2W_4=W_2W_4+W_6
$$
and thus (\ref{square}) since $W_2=\rho_2c_1=0$. Next let
$i:F\hookrightarrow K(\Z,4)\times K(\Z,6)$ denote the homotopy fibre of the induced map
$$
Sq^2\circ\rho_2+\rho_2:K(\Z,4)\times K(\Z,6)\to K(\Z_2,6).
$$
The relation (\ref{square}) implies that the map
$(c_2,c_3):BSU(\infty)\to K(\Z,4)\times K(\Z,6)$ is null--homotopic.
Consequently, $(c_2,c_3)$ lifts to a map $k:BSU(\infty)\to F$,
thereby inducing an isomorphism of homotopy groups
$\pi_i\big(BSU(\infty)\big)\to\pi_i(F)$ for $i\leq 7$ and a surjection for
$i=8$. By the exact homotopy sequence for fibrations we conclude
on one hand side that $\pi_4(F)=\Z$, $\pi_6(F)=2\Z$ and $\pi_i(F)=0$ for $i$
otherwise. On the other hand, the Chern class
$c_2:\pi_4\big(BSU(\infty)\big)\cong\widetilde{K}(S^4)=\Z\to \pi_4\big(K(\Z,4)\times K(\Z,6)\big)\cong H^4(S^4,\Z)=\Z$ is an isomorphism and
$c_3:\pi_4\big(BSU(\infty)\big)\cong\widetilde{K}(S^6)=\Z\to \pi_6\big(K(\Z,4)\times K(\Z,6)\big)\cong H^6(S^6,\Z)=\Z$ is multiplication
by~2. Since $M-B$ is 8--dimensional, it follows that the induced map $k_*:[(M-B)_+,BSU(\infty)]\to[(M-B)_+,F]$ is surjective. The horizontal row in the commutative diagram

\begin{center}
\begin{minipage}{11.5cm}
\unitlength0.4cm
\begin{picture}(2,5)
\put(11.15,0){$[(M-B)_+,BSU(\infty)]$}

\put(27,4.4){$H^6(M,\Z_2)$}

\put(21.8,5){$Sq^2\rho_2+\rho_2$}

\put(11.4,4.4){$H^4(M,\Z)\oplus H^6(M,\Z)$}

\put(-0.8,4.4){$[(M-B)_+,F]$}

\put(7.9,5){$i_*$}

\put(16.6,2.2){$(c_2,c_3)$}

\put(7.6,2.2){$k_*$}

\put(10.8,0.2){\vector(-2,1){7.5}}

\put(5.8,4.65){\vector(1,0){5.2}}

\put(21.6,4.65){\vector(1,0){5.2}}

\put(16.3,0.8){\vector(0,1){3.1}}
\end{picture}
\end{minipage}
\end{center}
is exact, hence ${\rm im}\,(c_2,c_3)={\rm im}\,i_*=\ker\,
(Sq^2\circ\rho_2+\rho_2)$.
\end{pf}

By assumption, $p_1\in H^4(M,\Z)$ is divisible by $6$ and therefore we can write $p_1=-6u$ for $u\in H^4(M,\Z)$. On the other hand,
$p_1=2q_1$, where $q_1$ is the first spin characteristic class which satisfies $\rho_2(q_1)=w_4$. Hence $Sq^2\rho_2(u)=Sq^2w_4=w_2w_4+w_6=0$,
and the previous lemma implies the existence of a stable complex vector bundle $\xi$ such that $c_1(\xi)=0$, $c_2(\xi)=u$ and
$c_3(\xi)=0$. From (\ref{pontrjagin}) it follows that
$p_1(\mf{su}(\xi))=p_1$, and since $w_2(\mf{su}(\xi))=0$,
$\mf{su}(\xi)$ and $TM$ are stably equivalent over the 4--skeleton
$M^{(4)}$~\cite{wo82b}. Then $\mf{su}(\xi)$ and $TM$ are stably equivalent over $M-B$ as the restriction map $KO(M-B)\to KO(M^{(4)})$ is injective. This follows from the exact sequence
$$
KO(M^{(i+1)},M^{(i)})\to KO(M^{(i+1)})\to KO(M^{(i)}).
$$
By definition
$KO(M^{(i+1)},M^{(i)})=\widetilde{KO}(M^{(i+1)}/M^{(i)})$ and $M^{(i+1)}/M^{(i)}$ is a disjoint union of spheres $S^{i+1}$. But
$\widetilde{KO}(S^{i+1})=0$ for $i=4,\,5$ and $6$ and therefore the map $KO(M^{(i+1)})\to KO(M^{(i)})$ is injective. Since $M=M^{(8)}$
is the disjoint union of $M^{(7)}$ and a finite number of open embedded discs, the assertion follows. Next we extend $\xi$ over $B$
to a stable bundle on $M$. The condition to be represented by a complex vector bundle $E$ of rank $3$ (which therefore is associated with a principal $SU(3)$--bundle) is $c_4(\xi)=0$. As pointed out
above, such a bundle exists if the spin index 
$$
\widehat{A}\cup{\rm ch}(\xi)[M]=3\widehat{A}[M]+p_1u/24+u^2/12
$$
is an integer, but this holds by assumption. Next $p_2(\mf{su}(\xi))=9u^2=p_2$ and as a consequence, $\mf{su}(\xi)$ is stably isomorphic to $TM$~\cite{wo82b}. Finally, two stably isomorphic oriented real vector bundles of rank 8 are isomorphic as $SO(8)$--bundles if they have the same Euler class. Since
$e\big(\mf{su}(\xi)\big)=0$, we conclude $TM\cong\mf{su}(\xi)$.
\end{pf}

\begin{cor}
If $M$ is closed and carries a $PSU(3)$--structure with vanishing
triality class, then $\widehat{A}[M]\in40\Z\mbox{ and }sgn(M)\in 640\Z$.
\end{cor}
%
%
%
%
%
\section{The twisted Dirac equation}
\label{Dirac}
In view of Hitchin's variational principle~\cite{hi01}, we adopt the following integrability condition, even if $M$ is not compact.

\begin{defn}\label{harmstruc}
A topological $PSU(3)$-- or $Sp(1)\cdot Sp(2)$--structure is called {\em harmonic}, if the defining $3$-- or $4$--form is closed and coclosed with respect to the metric it induces.
\end{defn}

\begin{rem}
The group $PSU(3)$ is the stabiliser of a totally symmetric $3$--tensor. See~{\rm \cite{gene64}} for an explicit description in terms of a $PSU(3)$--frame. Thus a symmetric $3$--tensor of the right algebraic type defines a topological $PSU(3)$--structure. Such structures, together with an integrability condition in some sense opposite to ours (cf. Remark~\ref{nearlyint}), were considered in~{\rm \cite{nu06}}.
\end{rem}

Our goal is to reformulate Definition~\ref{harmstruc} in terms of the supersymmetric maps associated with the topological $PSU(3)$-- and $Sp(1)\cdot Sp(2)$--structure. We show that the relevant $3$-- or $4$--form is harmonic if and only if the corresponding supersymmetric maps are in the kernel of the twisted Dirac operators $\Db_{\pm}:\Gamma(\Delta_{\pm}\otimes\Lambda^1)\to\Gamma(\Delta_{\mp}\otimes\Lambda^1)$.  Locally, these are given by
$$
\Db_{\pm}(\Psi\otimes X)=\sum
e_i\cdot\Psi\otimes\nabla_{e_i}X+e_i\cdot\nabla_{e_i}\Psi\otimes X,
$$
where $\nabla=\nabla^{LC}$ denotes the Levi--Civita connection as well as its lift to the spin bundle. 

To begin with, we recall the notion of intrinsic torsion. Consider an orbit in some $SO(n)$--representation space $V$ of the form $SO(n)/G$. A topological reduction of the principal frame bundle $P$ to a principal $G$--fibre bundle is characterised by a section $\gamma$ of the fibre bundle $P\times_{SO(n)}SO(n)/G$, of which we think as a section of the vector bundle $E=P\times_{SO(n)}V$. The Levi--Civita connection acts pointwise through $\mf{so}(n)\cong\Lambda^2$ on any section of $E$. In particular, since $\gamma$ is acted on trivially by its stabiliser algebra $\mf{g}$, 
\begin{equation}\label{torsion}
\nabla\gamma=T(\gamma),
\end{equation}
where $T$ is a section of the tensor bundle with fibre $\Lambda^1\otimes\mf{g}^{\perp}$, subsequently called the {\em torsion module}. The tensor field $T$ itself is referred to as the {\em intrinsic torsion} of the $G$--structure. If $\gamma$ is a $p$--form, this gives rise to the $G$--equivariant maps
$$
\mb{d}:X\otimes a\in\Lambda^1\otimes\mf{g}^{\perp} \mapsto X\wedge a(\gamma)\in\Lambda^{p+1},\,
\mb{d^*}:X\otimes a\in\Lambda^1\otimes\mf{g}^{\perp} \mapsto  X\llcorner a(\gamma)\in\Lambda^{p-1}.
$$
Since the differential operators $d$ and $d^*$ are induced by skew--symmetrisation and minus the contraction of the Levi--Civita connection, we deduce from~(\ref{torsion}) that $\mb{d}(T)=d\gamma$ and $\mb{d^*}(T)=-d^*\gamma$. Consequently, $\gamma$ is harmonic if and only if $T\in\ker\mb{d}\cap\ker\mb{d^*}$. On the other hand, for $\gamma$ a spinor--valued $1$--form, we can consider the equivariant map
$$
\Dbb:\,\,X\otimes a\in\Lambda^1\otimes\mf{g}^{\perp}\mapsto\mu\big(a(\gamma)\otimes X\big)\in\Delta\otimes\Lambda^1.
$$
Here, $a$ acts via the induced action of $\mf{so}(n)$ on $\Delta_{\pm}\otimes\Lambda^1$, i.e.
$$
X\wedge Y(\Psi_{\pm}\otimes Z)=\frac{1}{4}(X\cdot Y-Y\cdot X)\cdot\Psi_{\pm}\otimes Z+\Psi_{\pm}\otimes\big(X(Z)Y-Y(Z)X\big),
$$ 
and Clifford multiplication takes $\Psi_{\pm}\otimes Z\otimes X\in\Delta_{\pm}\otimes\Lambda^1\otimes\Lambda^1$ to $X\cdot\Psi_{\pm}\otimes Z\in\Delta_{\mp}\otimes\Lambda^1$. Hence
$\Dbb(T)=\Db\gamma$, and our task consists in showing that $\ker\Dbb=\ker\mb{d}\cap\ker\mb{d^*}$. 

By equivariance, the kernels of $\Dbb$, $\mb{d}$ and $\mb{d^*}$ can be computed using Schur's Lemma and $G$--representation theory. From a technical point of view, the $Sp(1)\cdot Sp(2)$--case is a lot easier to deal with, so we start with this one. Here, the invariant 4--form $\Omega$ is self--dual, so we only need to show $\ker\Dbb=\ker\mb{d}$. First, we decompose the torsion module $\Lambda^1\otimes\big(\mf{sp}(1)\oplus\mf{sp}(2)\big)^{\perp}=[1,0,1]\otimes[1,1,1]$ into $Sp(1)\cdot Sp(2)$--irreducibles which yields
$$
\Lambda^1\otimes\big(\mf{sp}(1)\oplus\mf{sp}(2)\big)^{\perp}=[{\textstyle\frac{1}{2}},{\textstyle\frac{1}{2}},1]\oplus[{\textstyle\frac{1}{2}},{\textstyle\frac{3}{2}},2]\oplus[{\textstyle\frac{3}{2}},{\textstyle\frac{1}{2}},1]\oplus[{\textstyle\frac{3}{2}},{\textstyle\frac{3}{2}},2].
$$
On the other hand, we find for the target spaces of $\Dbb$ and $\mb{d}$ (cf. Proposition~\ref{sp(1)sp(2)decomp})
$$
\Delta_-\otimes\Lambda^1=2
[{\textstyle\frac{1}{2}},{\textstyle\frac{1}{2}},1]\oplus[{\textstyle\frac{1}{2}},{\textstyle\frac{3}{2}},2]\oplus[{\textstyle\frac{3}{2}},{\textstyle\frac{1}{2}},1],\quad\Lambda^5=[{\textstyle\frac{1}{2}},{\textstyle\frac{1}{2}},1]\oplus[{\textstyle\frac{1}{2}},{\textstyle\frac{3}{2}},2]\oplus[{\textstyle\frac{3}{2}},{\textstyle\frac{1}{2}},1].
$$

\begin{thm}\label{sp1sp2thm}
The topological $Sp(1)\cdot Sp(2)$--structure $(M^8,\Omega)$ is harmonic if one of the following equivalent statements holds.

{\rm (i)} $d\Omega=0$.

{\rm (ii)} If $\sigma_+\in\Gamma(\Delta_+\otimes\Lambda^1)$ is the corresponding supersymmetric map, then $\Db_+(\sigma_+)=0$.

{\rm (iii)} The intrinsic torsion $T$ takes values in $[{\textstyle\frac{3}{2}},{\textstyle\frac{3}{2}},2]$. 
\end{thm}

\begin{pf}
By Schur's Lemma, it is enough to evaluate the maps $\Dbb$ and $\mb{d}$ for a sample vector of a given module to check whether or not it maps non--trivially. The operator $\mb{d}$ is known to be surjective~\cite{sa89} whence $\ker\mb{d}=[{\textstyle\frac{3}{2}},{\textstyle\frac{3}{2}},2]$. On the other hand, $[{\textstyle\frac{3}{2}},{\textstyle\frac{3}{2}},2]\subset\ker\Dbb$ and we are left with showing that the remaining irreducible modules in $\Lambda^1\otimes\big(\mf{sp}(1)\oplus\mf{sp}(2)\big)^{\perp}$ map non--trivially. To that end we consider the operator $L:\Lambda^3\to\Lambda^3$ built out of $\Dbb$, compounded with $id\otimes\sigma_+$ and the projection $\Delta_-\otimes\Delta_+\to\Lambda^3$ given by $\sum q(\Psi_-,e_I\cdot\Psi_+)e_I$ with respect to some $Sp(1)\cdot Sp(2)$--frame. Since all the representations involved are of real type, $L$ restricted to an irreducible module of $\Lambda^3$ is multiplication by a real scalar, possibly zero. To begin with, we map the element $4e_1\llcorner\Omega\in[{\textstyle\frac{1}{2}},{\textstyle\frac{1}{2}},1]$ into $\Lambda^1\otimes\Lambda^2$ in the natural way, namely $x\wedge y\wedge z\mapsto(x\otimes y\wedge z+cyc.\,\,perm.)/3$. 
Projecting the second factor onto $[1,1,1]=\big(\mf{sp}(1)\oplus\mf{sp}(2)\big)^{\perp}$ by means of the projection operator given in Proposition~\ref{2formssp1sp2} yields the element
\begin{eqnarray*}
t_1 & = & e_2\otimes (e_{12}-e_{34}-e_{56}+e_{78})+e_3\otimes(e_{13}+e_{24}-e_{57}-e_{68})+\\
& & e_4\otimes(e_{14}-e_{23}-e_{58}+e_{67})+e_5\otimes(3e_{15}-e_{26}-e_{37}-e_{48})+\\
& & e_6\otimes(3e_{16}+e_{25}-e_{38}+e_{47})+e_7\otimes(3e_{17}+e_{28}+e_{35}-e_{46})+\\
& & e_8\otimes(3e_{18}-e_{27}+e_{36}+e_{45}) \in[{\textstyle\frac{1}{2}},{\textstyle\frac{1}{2}},1]\subset\Lambda^1\otimes[1,1,1]. 
\end{eqnarray*}
Using the matrix representation of $\sigma_+$ in Appendix~\ref{susymatrix}, and evaluating $L$ on $t_1$ yields
$$
L(t_1)=-6e_{234}+2e_{256}-2e_{278}+2e_{357}+2e_{368}+2e_{458}-2e_{467}=2t_1.
$$
In particular, $t_1$ maps non--trivially under $\Dbb$. Next we turn to the module $[{\textstyle\frac{3}{2}},{\textstyle\frac{1}{2}},1]$. 
Using the weight vectors provided in~(\ref{weightssp1sp2}), the vector 
\begin{eqnarray*}
x_{(\alpha+\beta_1)/2+\beta}\wedge x_{(\alpha-\beta_1)/2}\wedge x_{(\alpha+\beta_1)/2} & = & e_{157}+ie_{158}-ie_{167}+e_{168}-\\
& & ie_{257}+e_{258}-e_{267}-ie_{268}
\end{eqnarray*}
is of the highest weight occurring in $\Lambda^3\otimes\C$. Hence it is actually the weight vector of 
$({\textstyle\frac{3}{2}},{\textstyle\frac{1}{2}},1)\subset \Lambda^3\otimes\C$. Proceeding as before yields the sample vector $t_2\in({\textstyle\frac{3}{2}},{\textstyle\frac{1}{2}},1)\subset\Lambda^1\otimes[1,1,1]\otimes\C$ which is mapped to $20t_2$ under $L$. For the remaining module $[{\textstyle\frac{1}{2}},{\textstyle\frac{3}{2}},2]$ we consider the vector $-e_{134}/3+e_{178}$. Wedging with $\Omega$ yields zero, so it is necessarily contained in $[{\textstyle\frac{1}{2}},{\textstyle\frac{3}{2}},2]\oplus[{\textstyle\frac{3}{2}},{\textstyle\frac{1}{2}},1]$. Mapping it into $\Lambda^1\otimes[1,1,1]$ gives the vector $t_3$ whose projections on $[{\textstyle\frac{1}{2}},{\textstyle\frac{3}{2}},2]$ and $[{\textstyle\frac{3}{2}},{\textstyle\frac{1}{2}},1]$ we denote by $t_{3a}$ and $t_{3b}$. The image under $L$ is
\begin{eqnarray*}
L(t_{3a}+t_{3b}) & = & c\cdot t_{3a}+20\cdot t_{3b}\\
& = & (-12e_{134}-8e_{156}+44e_{178}+4e_{358}-4e_{367}-4e_{457}-4e_{468})/3 
\end{eqnarray*}
Since $L(t_3)-20t_3\not=0$, we see that $[{\textstyle\frac{1}{2}},{\textstyle\frac{3}{2}},2]$ also maps non--trivially (applying $L$ again to $L(t_3)-20t_3$ shows that restricted to this module, $L$ is actually multiplication by $12$). Hence $\ker\Dbb=[{\textstyle\frac{3}{2}},{\textstyle\frac{3}{2}},2]$ which proves the theorem.
\end{pf}

\begin{cor}
Let $\sigma=\C^2$ be the vector representation of $Sp(1)=SU(2)$ and $\lambda^r_s$ the class of irreducible $Sp(2)$--representations introduced after Proposition~\ref{2formssp1sp2}. If $L:\Lambda^3\to\Lambda^3$ denotes the $Sp(1)\cdot Sp(2)$--invariant map defined in the proof of Theorem~\ref{sp1sp2thm}, the irreducible $Sp(1)\cdot Sp(2)$--modules can be characterised as follows:
$$
\begin{array}{lclcl}
\,[\sigma\otimes\lambda^1_0] & = & [\frac{1}{2},\frac{1}{2},1] & = &\{\alpha\in\Lambda^3\,|\,L(\alpha)=2\alpha\}\\
\,[\sigma\otimes\lambda^3_1] & = & [\frac{1}{2},\frac{3}{2},2] & = &\{\alpha\in\Lambda^3\,|\,L(\alpha)=12\alpha\}\\
\,[\sigma^3\otimes\lambda^1_0] & = & [\frac{3}{2},\frac{1}{2},1] & = &\{\alpha\in\Lambda^3\,|\,L(\alpha)=20\alpha\}
\end{array}
$$
\end{cor}

Next we turn to $PSU(3)$. The situation here is more involved not only because the defining form $\rho$ is not self--dual anymore, but also due to the presence of modules with multiplicities greater than one. 

Again we begin by decomposing the torsion module. Let $\wedge:\Lambda^1\otimes\Lambda^2_{20}\to\Lambda^3$ denote the natural skewing map. Then $\Lambda^1\otimes\Lambda^2_{20}\cong\ker\wedge\oplus\rho^{\perp}$,
where $\rho^{\perp}=[1,1]\oplus\llbracket1,2\rrbracket\oplus[2,2]$ is the orthogonal complement of $\rho$ in $\Lambda^3$. Moreover, the natural contraction map $\llcorner:\ker\wedge\subset\Lambda^1\otimes\Lambda^2_{20}\to\Lambda^1$ splits $\ker\wedge$ into a direct sum isomorphic to $\ker\,\llcorner\oplus\Lambda^1$, where $\ker\,\llcorner\cong[2,2]\oplus\llbracket2,3\rrbracket$. Consequently, the complexification of $\Lambda^1\otimes\Lambda^2_{20}$ is the direct sum of 
$$
\begin{array}{ccccl}
& &\Lambda^1\otimes\Lambda^2_{10+} & = & (1,1)_+\oplus(1,2)\oplus(2,2)_+\oplus(2,3)\\
\Lambda^1\otimes\Lambda^2_{20}\otimes\C & = & \mbox{\hspace{-14pt}}\oplus & & \\
& &\Lambda^1\otimes\Lambda^2_{10-} & = & (1,1)_-\oplus(2,1)\oplus(2,2)_-\oplus(3,2).
\end{array}
$$
The modules $(1,1)_{\pm}$ and $(2,2)_{\pm}$ have non--trivial projections to both $\ker\wedge$ and $\rho^{\perp}$. In particular, they map non--trivially under $\wedge$. With the decomposition of the target spaces of $\Dbb_{\pm}=\Dbb_{|\Delta_{\pm}\otimes\Lambda^1}$, $\mb{d}$ and $\mb{d^*}$, namely
$$
\Delta_{\mp}\otimes\Lambda^1=\mb{1}\oplus2[1,1]\oplus\llbracket1,2\rrbracket\oplus[2,2],\quad\Lambda^4=2[1,1]\oplus2[2,2],\quad\Lambda^2=[1,1]\oplus\llbracket1,2\rrbracket,
$$
we can now prove the analogue of Theorem~\ref{sp1sp2thm}.

\begin{thm}\label{psu3thm}
The topological $PSU(3)$--structure $(M^8,\rho)$ is harmonic if one of the following equivalent statements holds.

{\rm (i)} $d\rho=0$ and $d\star_{\rho}\rho=0$.
	
{\rm (ii)} If  $\sigma_{\pm}\in\Gamma(\Delta_{\pm}\otimes\Lambda^1)$ are the corresponding supersymmetric maps, then $\Db_{\pm}(\sigma_{\pm})=0$.
	
{\rm (iii)} The intrinsic torsion $T$ takes values in $\llbracket2,3\rrbracket$. 
\end{thm}

\begin{rem}\label{proof}
The implication $(i)\Rightarrow (ii)$ was already asserted in~{\rm\cite{hi01}}. However, the proof is inconclusive. Firstly, some of the sample vectors provided in the proof are not contained in the right module. For instance, $x_{\alpha_1}\otimes x_{\alpha_1}\wedge x_{\alpha_2}$ is not contained in $(2,1)\subset\Lambda^1\otimes\Lambda^2_{10-}$ as claimed: Using the author's notation, applying $[x_{\alpha_1},\cdot]$ yields $x_{\alpha_1}\otimes x_{\alpha_1}\wedge x_{\alpha_1+\alpha_2}$ from general properties of root vectors. Moreover, due to the presence of modules with multiplicity two, the modules $(1,1)_{\pm}$ and $(2,2)_{\pm}$ can map non--trivially under say $\mb{d}$ while $\ker\mb{d}$ still contains a component isomorphic to $(1,1)$ or $(2,2)$.
\end{rem}

\begin{pf}
We first establish the equivalence between (i) and (iii) and start by determining the kernel of $\mb{d^*}$. By Proposition~\ref{liealgsu3}, $a_{\pm}(\rho)=\pm\sqrt{3}i\star(a_{\pm}\wedge\rho)$ for any $a_{\pm}\in\Lambda^2_{10\pm}$. It follows from  complexifying that restricted to the $PSU(3)$--invariant modules $\Lambda^1\otimes\Lambda^2_{10\pm}$,
\begin{equation}\label{dstar_van}
\mb{d^*}(\sum e_j\otimes a_j^{\pm})=\mp\sqrt{3}i\sum e_j\llcorner\star(a^{\pm}_j\wedge\rho)=\pm\sqrt{3}i\sum \star(e_j\wedge a^{\pm}_j\wedge\rho).
\end{equation}
In virtue of the remarks above, the kernel of the skewing map $\Lambda^1\otimes\Lambda^2_{\pm}$ is isomorphic to $(2,3)$ and $(3,2)$. Hence,~(\ref{dstar_van}) vanishes if and only if $\sum e_i\wedge a_i^{\pm}$ lies in $\mb{1}\oplus[2,2]$, the kernel of the map which wedges 3--forms with $\rho$. Invoking Schur's Lemma, $\ker\mb{d^*}\cong[1,1]\oplus2[2,2]\oplus\llbracket2,3\rrbracket$, 
where the precise embedding of $[1,1]$ will be of no importance to us.

Next we consider the operator $\mb{d}$. If we can show that it is surjective, then $\ker\mb{d}\cong\llbracket1,2\rrbracket\oplus\llbracket2,3\rrbracket$ and consequently, the kernels of $\mb{d^*}$ and $\mb{d}$ intersect in $\llbracket2,3\rrbracket$. Let $\iota_{\rho^{\perp}}$ denote the injection of $\rho^{\perp}$ into $\Lambda^1\otimes\Lambda^2_{20}$ obtained by projecting the natural embedding of $\Lambda^3$ into $\Lambda^1\otimes\Lambda^2$. We first prove the relation (cf. Section~\ref{triality} for the definition of $c_q$)
\begin{equation}\label{surjd}
c_3(\alpha)=\frac{1}{2}\mb{d}\big(\iota_{\rho^{\perp}}(\alpha)\big),\quad \alpha\in\rho^{\perp}\subset\Lambda^3
\end{equation}
which shows that $\ker c_4\subset\im\mb{d}$. By~(\ref{b_decomp}), the kernel of $c_3$ is isomorphic to $\mb{1}\oplus\llbracket1,2\rrbracket$, so the claim needs only to be verified for the module $[1,1]\oplus[2,2]$ in $\Lambda^3$. A sample vector is obtained by 
\begin{equation}\label{sample}
p_3(e_{128})=\alpha_8\oplus\alpha_{27}={\textstyle\frac{1}{8}}(5e_{128}+\sqrt{3}e_{345}+\sqrt{3}e_{367}-2e_{458}+2e_{678}),
\end{equation}
where $p_3=c^*_4c_3$. That both components $\alpha_8$ and $\alpha_{27}$ are non--trivial can be seen as follows. Restricting $p_3$ to $\Lambda^3_8$ and $\Lambda^3_{27}$ is multiplication by real scalars $x_1$ and $x_2$ since the modules are representations of real type. If one, say $x_1$, were to vanish, then $p^2_3(e_{128})=p_3(\alpha_{27})=x_2\cdot\alpha_{27}$. However
$$
p^2_3(e_{128})={\textstyle\frac{1}{64}}(39e_{128}+7\sqrt{3}e_{345}+7\sqrt{3}e_{367}-18e_{458}+18e_{678})
$$
which is not a multiple of~(\ref{sample}). Moreover, we have indeed 
\begin{eqnarray*}
c_3p_3(e_{128}) & = & {\textstyle\frac{1}{32}}(7\sqrt{3}e_{1245}+7\sqrt{3}e_{1267}-9e_{1468}-9e_{1578}+9e_{2478}-9e_{2568})\\
& = & {\textstyle\frac{1}{2}}\mb{d}\big(\iota_{\rho^{\perp}}p_3(e_{128})\big)
\end{eqnarray*} 
which proves~(\ref{surjd}). For the inclusion $\im c^*_5\subset\im\mb{d}$ we consider the vector $e_1\otimes e_{18}$ in $\ker\wedge$. Then $\mb{d}(e_1\otimes e_{18})=-e_{1238}/2-e_{1478}/4+e_{1568}/4$ takes values in both components of $\im c^*_5\subset\Lambda^4$, since
$c^*_4\mb{d}(e_1\otimes e_{18})=0$ and otherwise 
$$
c^*_5c_4\mb{d}(e_1\otimes e_{18})={\textstyle\frac{1}{32}}(-10e_{1238}-5e_{1478}+5e_{1568}+3e_{2468}+3e_{2578}+3e_{3458}-3e_{3678})
$$ 
would be a multiple of $\mb{d}(e_1\otimes e_{18})$. Hence $\mb{d}$ is surjective and the equivalence between (i) and (iii) is established.

Finally, we turn to the twisted Dirac equation. We have to prove that
$$
\ker\Dbb_+\cap\ker\Dbb_-=\llbracket2,3\rrbracket=\ker\mb{d}\cap\ker \mb{d^*}.
$$
The intersection $\ker\Dbb_+\cap\ker\Dbb_-$ contains at least the module $\llbracket2,3\rrbracket$. First we show that $\llbracket1,2\rrbracket$ is not contained in this intersection by taking the vector 
\begin{eqnarray*}
\tau_{\llbracket1,2\rrbracket} & = & id\otimes\pi^2_{20}\big(\iota_{\rho^{\perp}}c_2(4e_{18})\big)\\
& = & -\sqrt{3}e_1\otimes e_{45}-\sqrt{3}e_1\otimes e_{67}+2e_2\otimes e_{38}-2e_3\otimes e_{28}+
\sqrt{3}e_4\otimes e_{15}\\
& &+e_4\otimes e_{78}-\sqrt{3}e_5\otimes e_{14}-e_5\otimes e_{68}+\sqrt{3}e_6\otimes e_{17}+e_6\otimes e_{58}\\
& &-\sqrt{3}e_7\otimes e_{16}-e_7\otimes e_{48}+2e_8\otimes e_{23}+e_8\otimes e_{47}-e_8\otimes e_{56}.
\end{eqnarray*}
A straightforward, if tedious, computation shows $\Dbb_{\pm}(\tau_{\llbracket1,2\rrbracket})\not=0$. For the remainder of the proof, it will again be convenient to complexify the torsion module $\Lambda^1\otimes\Lambda^2_{20}$ and to consider $(1,1)_{\pm}$ and $(2,2)_{\pm}$. The invariant 3--form $\rho$ induces equivariant maps $\rho_{\mp}:\Delta_{\pm}\to\Delta_{\mp}$ whose matrices with respect to the choices made in~(\ref{spin8rep}) are given by~(\ref{rhomap}) for $\rho_+$, and by its transpose for $\rho_-$. Schur's Lemma implies
\begin{equation}\label{constant}
\Dbb_-\big((2,2)_+\big)=z\cdot\rho_-\otimes id\circ\Dbb_+\big((2,2)_+\big)
\end{equation}
for a complex scalar $z$. Since the operators $\Dbb_{\pm}$ are real and $(2,2)_-$ is the complex conjugate of $(2,2)_+$, the same relation holds for $(2,2)_-$ with $\bar{z}$. The vector $\tau_0=6(e_1\otimes e_{18}-e_2\otimes e_{28})$ is clearly in $\ker\,\llcorner\subset\ker\wedge$ and projecting the second factor to $\Lambda^2_{10+}$ yields 
\begin{eqnarray*}
id\otimes\pi^2_{10+}(\tau_0) & = & e_1\otimes(3e_{18}+i\sqrt{3}e_{23}-i\sqrt{3}e_{47}+i\sqrt{3}e_{56})+\\
& & e_2\otimes(i\sqrt{3}e_{13}-3e_{28}+i\sqrt{3}e_{46}+i\sqrt{3}e_{57}).
\end{eqnarray*}
Since any possible component in $(2,3)$ gets killed under $\Dbb_{\pm}$, we can plug this into~(\ref{constant}) to find $z=(1+i\sqrt{3})/8$. This shows that $(2,2)_{\pm}$ maps non--trivially under $\Dbb$. On dimensional grounds, $\ker\Dbb_{\pm}$ therefore contains the module $(2,2)$ with multiplicity one. Their intersection, however, is trivial, for suppose otherwise. Let $(2,2)_0$ denote the corresponding copy in $\ker\Dbb_+$.
It is the graph of an isomorphism $P:(2,2)_+\to(2,2)_-$ since it intersects $(2,2)_{\pm}$ trivially. Now if $\tau=\tau_+\oplus P\tau_+\in(2,2)_0$ were in $\ker\Dbb_-$, then
\begin{eqnarray*}
\Dbb_-(\tau_+\oplus P\tau_+) & = & z\cdot\rho\otimes id\circ\Dbb_+(\tau_+)\oplus \bar{z}\cdot\rho\otimes id\circ\Dbb_+(P\tau_+)\\
& = & \rho\otimes id\circ\Dbb_+(z\cdot \tau_+\oplus \bar{z}\cdot P\tau_+)\\
& = & 0.
\end{eqnarray*}
Consequently, $z\cdot \tau_+\oplus \bar{z}\cdot P\tau_+\in\ker\Dbb_+$, that is, $\bar{z}\cdot P\tau_+=Pz\cdot \tau_+$ or $\bar{z}=z$ which is a contradiction. This shows that (a) the kernels of $\Dbb_{\pm}$ intersect at most in $2(1,1)\oplus\llbracket2,3\rrbracket$ and (b)  the conditions $\Db_+(\sigma_+)=0$ or $\Db_-(\sigma_-)=0$ on their own are not sufficient to guarantee the close-- and cocloseness of $\rho$. The same argument also applies to $(1,1)_{\pm}$. However, since $(1,1)$ appears twice in $\Delta_{\pm}\otimes\Lambda$, we first need to project onto $\Delta_{\mp}\cong(1,1)$ via Clifford multiplication before asserting the existence of a complex scalar $z$ such that
$$
\mu_+\circ\Dbb_-\big((1,1)_+\big)=z\cdot\rho_+\big(\mu_-\circ\Dbb_+\big((1,1)_+\big)\big).
$$
For the computation of $z$, we can use the vector 
$$
2\sqrt{3}ie_1\otimes\pi^2_{10+}(e_{18})= e_1\otimes(\sqrt{3}ie_{18}-e_{23}+e_{47}- e_{56})\in(1,1)_+\oplus(2,2)_+\oplus(2,3),
$$
as possible non--trivial components in $(2,2)_+\oplus(2,3)$ get killed under $\mu_{\mp}$. We find $z=2(1-\sqrt{3}i)$ which as above shows that $(1,1)$ occurs with multiplicity at most one in $\ker\Dbb_{\pm}$, and that it is not contained in their intersection. Consequently, $\ker\Dbb_+\cap\ker\Dbb_-=\llbracket2,3\rrbracket$, which proves the equivalence between (ii) and (iii).
\end{pf}

\begin{rem}\label{nearlyint}
{\rm (i)} From (iii) it follows that our harmonicity condition on a topological $PSU(3)$--structure can be seen as the extreme opposite of Nurowski's notion of restricted nearly integrable $PSU(3)$--structures, where $T$ takes values in the remaining modules isomorphic to $[1,1]$, $\llbracket1,2\rrbracket$ and $[2,2]$~{\rm\cite{nu06}}.

{\rm (ii)} According to the decomposition $\Delta_{\pm}\otimes\Lambda^1=\Delta_{\mp}\oplus\ker\mu_{\pm}$, the twisted Dirac operators $\Db_{\pm}:\Gamma(\Delta_{\pm}\otimes\Lambda^1)\to\Gamma(\Delta_{\mp}\otimes\Lambda^1)$ take the shape~{\rm\cite{wa91}}
\begin{equation}\label{dodecomp}
\Db_{\pm}=\left(\begin{array}{cc} -\frac{3}{4}\iota_{\pm}\circ D_{\mp}\circ \iota^{-1}_{\mp} &
2\iota_{\pm}\circ\db_{\pm}\\\frac{1}{4}P_{\mp}\circ \iota^{-1}_{\mp} & Q_{\pm}\end{array}\right).
\end{equation}
Here, $\iota_{\pm}:\Gamma(\Delta_{\pm})\to \Gamma(\Delta_{\mp}\otimes\Lambda^1)$ is the embedding given by $\iota_{\pm}(\Psi_{\pm})(X)=-X\cdot\Psi_{\pm}/8$, and $D_{\pm}:\Gamma(\Delta_{\pm})\to\Gamma(\Delta_{\mp})$ and $P_{\pm}:\Gamma(\Delta_{\pm})\to \ker\mu_{\pm}$ denote the usual Dirac-- and Penrose--operator. Further, $\db_{\pm}:\Delta_{\pm}\otimes\Lambda^1\to\Delta_{\pm}$ is the twisted co--differential and $Q:\ker\mu\to\ker\mu$ the Rarita--Schwinger operator. In particular, we see that the supersymmetric maps of a harmonic $PSU(3)$-- or $Sp(1)\cdot Sp(2)$--structure define {\em Rarita--Schwinger fields}, spin $3/2$ particles satisfying the relativistic field equation $Q\sigma=0$ in physicists' language~{\rm\cite{rasc41}}. Further, if $\Psi\in\Gamma(\Delta_+)\subset\Gamma(\Delta_-\otimes\Lambda^1)$ denotes the invariant spinor coming from a topological $Spin(7)$--structure $(M^8,\Omega)$ with $Spin(7)$--invariant $4$--form $\Omega$, then $\Db_-(\Psi)=0$ if and only if $D_+(\Psi)=0$ and $P_+(\Psi)=0$, hence if and only if $\nabla\Psi=0$. This implies that the holonomy of the induced metric is contained in $Spin(7)$, which by~{\rm\cite{br87}} is equivalent to $d\Omega=0$. In this way, $Sp(1)\cdot Sp(2)$--structures appear on an equal footing with $Spin(7)$--manifolds. However, in contrast to these, closeness of the $Sp(1)\cdot Sp(2)$--invariant $4$--form does not imply the holonomy to be contained in $Sp(1)\cdot Sp(2)$ (cf. Salamon's counterexample in~{\rm\cite{sa01}} given in the next section), although this is true for $Sp(1)\cdot Sp(k)$--structures on $M^{4k}$ with $k\geq3$~{\rm\cite{sw89}}.
\end{rem}

Next we derive integrability conditions for harmonic $PSU(3)$-- and $Sp(1)\cdot Sp(2)$--structures on the Ricci tensor: According to Proposition 2.8 in~\cite{wa91}, 
$$
(D\circ\db-\db\circ \Db)(\gamma)=\frac{1}{2}p(\gamma\circ Ric)
$$
for any $\gamma=\sum_ie_i\otimes\gamma_i\in\Gamma(T^*M\otimes\Delta)$, where $\gamma\circ Ric=\sum_{i,j}Ric_{ij}e_i\otimes\gamma_j$.
Integrability implies
$$
p(\sigma\circ Ric)=\sum_{i,j}Ric_{ij}e_i\cdot\sigma_j=0.
$$
Hence, in our case $Ric$ is in the kernel of the map
$$
A\in\odot^2\mapsto
p(A\circ\sigma)=\sum_{i,j}A_{ij}e_i\sigma_j\in\Delta,
$$
which is invariant under the stabiliser of $\sigma$. For a $PSU(3)$--structure, we have $\odot^2=\mathbf{1}\oplus [1,1]\oplus
[2,2]$ and $\Delta_{\pm}\cong[1,1]$, so $Ric$ vanishes on the module $[1,1]$, a fact previously noted in~\cite{hi01}. For an $Sp(1)\cdot Sp(2)$--structure, $\sigma=\sigma_+\in\Delta_+\otimes\Lambda$, hence $\odot^2=\mathbf{1}\oplus [0,1,1]\oplus
[1,1,2]$ and $\Delta_-\cong[2,0,0]\oplus[0,1,1]$. Since this map is non--trivial, $Ric$ vanishes on the module $[0,1,1]$.

\begin{prop}\label{riccitensor}\hfill\newline
{\rm (i)} If $g$ is a metric induced by a harmonic $PSU(3)$--structure, then $Ric$ vanishes on the ($8$--dimensional) component $[1,1]$.

{\rm (ii)} If $g$ is a metric induced by a harmonic $Sp(1)\cdot Sp(2)$--structure, then $Ric$ vanishes on the ($5$--dimensional) component $[0,1,1]$.
\end{prop}

\begin{rem}
As shown by the examples in the next section there exist (compact) harmonic $PSU(3)$-- or $Sp(1)\cdot Sp(2)$--structures whose Ricci tensor takes values in the remaining modules (i.e. $\mathbf{1}$ and $[2,2]$ for $PSU(3)$ and $\mathbf{1}$ and $[1,1,2]$ for $Sp(1)\cdot Sp(2)$). In particular, unlike harmonic $G_2$-- or $Spin(7)$--structures which are Ricci--flat, harmonic $PSU(3)$-- or $Sp(1)\cdot Sp(2)$--structures are not necessarily Einstein (and a fortiori, not Ricci--flat).
\end{rem}
%
%
%
%
%
\section{Examples}
\label{examples}
Let $B$ denote the Killing form of $\mf{su}(3)$. Then $\rho(X,Y,Z)=-4B([X,Y],Z)/3$ provides $\mf{su}(3)$,
the symmetric space $SU(3)=SU(3)\times SU(3)/SU(3)$ and its
non--compact dual $SL(3,\C)/SU(3)$ with topological $PSU(3)$--structures. Since $\rho$ and $\star\rho$ are
bi--invariant, they are closed and induce harmonic $PSU(3)$--metrics which are Einstein and respectively of zero, positive and negative scalar curvature. In fact, by standard symmetric space theory, the holonomy is contained in $PSU(3)$, whence $\nabla\rho=0$. Conversely, any irreducible $PSU(3)$--structure with parallel $3$--form is either locally symmetric or flat as follows from inspection of Berger's list. In this section, we construct examples with non--parallel $3$--form.

\medskip
{\bf Local examples.} The first example is built out of a hyperk\"ahler 4--manifold $M^4$ with a triholomorphic vector field. Let $\theta$ be a 1--form on $\R^3$ and $U\equiv
U(x,y,z)$ be a strictly positive harmonic function on some domain
$D\subset\R^3$ with $dU=\star d\theta$. By the Gibbons--Hawking ansatz~\cite{ar02},~\cite{giha78}, the metric 
\begin{equation}\label{hypmet}
h=U(dx^2+dy^2+dz^2)+\frac{1}{U}(dt+\theta)^2
\end{equation}
on $D\times\R$ is hyperk\"ahler with associated K\"ahler forms given by
\begin{eqnarray*}
\omega_1^- & = & Udy\wedge dz + dx\wedge (dt+\theta)\\
\omega_2^- & = & Udx\wedge dy + dz\wedge (dt+\theta)\\
\omega_3^- & = & Udx\wedge dz - dy\wedge (dt+\theta).
\end{eqnarray*}
The vector field $X=\frac{\partial}{\partial t}$ is triholomorphic, hence induces an infinitesimal transformation which preserves
each of the three complex structures associated with $\omega_1^-$,
$\omega_2^-$ or $\omega_3^-$. Conversely, a hyperk\"ahler metric on a
4--dimensional manifold which admits a triholomorphic vector field is
locally of the form~(\ref{hypmet}).

Let us now define the 2--form $\omega_3^+$
by changing the sign in $\omega_3^-$, that is
$$ 
\omega_3^+=Udx\wedge dz + dy\wedge(dt+\theta).
$$ 
This 2--form is closed if and only if $U\equiv U(x,z)$, for $d\omega_3^+=0$ implies
\begin{equation}\label{de45}
d(Udx\wedge dz)=d\big(dy\wedge(dt+\theta)\big),
\end{equation}
whence
$$
d\omega_3^+=2d(Udx\wedge dz)=2\frac{\partial U}{\partial y}dy\wedge dx\wedge dz.
$$
Pick such a $U$ and take the standard coordinates $x_1,\ldots,x_4$
of the Euclidean space $(\R^4,g_0)$. Put
$$
\begin{array}{llll} e^1=dx_1, & e^2=dx_2, & e^3=dx_3, & e^8=dx_4, \\
e^4=\sqrt{U}dy, & e^5=-\frac{1}{\sqrt{U}}(dt+\theta), &
e^6=-\sqrt{U}dx, & e^7=\sqrt{U}dz,
\end{array}
$$
which we take as an orthonormal coframe on $M^4\times \R^4$. As before, we shall drop any distinction between vector fields and $1$--forms in the presence of a metric. Endowed
with the orientation defined by $(e_4,\ldots,e_7)$, the forms
$\omega_i^-$ are anti--self--dual on $M^4$, while the forms
$\omega_1^+=Udy\wedge dz - dx\wedge (dt+\theta)$,
$\omega_2^+=Udx\wedge dy - dz\wedge (dt+\theta)$ and
$\omega_3^+$ are self--dual, so
\begin{equation}\label{wedge}
\omega_i^{\pm}\wedge\omega_j^{\mp}=0,\quad
\omega_i^{\pm}\wedge\omega_j^{\pm}=\pm2\delta_{ij}e_{4567}.
\end{equation}
The 3--form
\begin{equation}\label{rho1s}
\rho = \frac{1}{2}e_{123}+\frac{1}{4}e_1\wedge\omega_1^-+\frac{1}{4}e_2\wedge\omega_2^-+\frac{1}{4}e_3\wedge\omega_3^-+\frac{\sqrt{3}}{4}e_8\wedge\omega_3^+
\end{equation}
defines a $PSU(3)$--structure which is closed by design. Moreover, the same holds for
\begin{eqnarray}
\star\rho & = & \phantom{+}\frac{1}{2}e_{45678}-\frac{1}{4}\omega_1^-\wedge
e_{238}+\frac{1}{4}\omega_2^-\wedge e_{138}-\frac{1}{4}\omega_3^-\wedge
e_{128}+\frac{\sqrt{3}}{4}\omega_3^+\wedge e_{123}\nonumber\\
& = & \phantom{+}\frac{1}{2}Udx\wedge dz\wedge
dy\wedge(dt+\theta)\wedge dx_4-\frac{1}{4}\omega_1^-\wedge dx_2\wedge
dx_3\wedge dx_4\nonumber\\ & & +\frac{1}{4}\omega_2^-\wedge dx_1\wedge
dx_3\wedge dx_4-\frac{1}{4}\omega_3^-\wedge dx_1\wedge dx_2\wedge
dx_4\nonumber\\
& &+\frac{\sqrt{3}}{4}\omega_3^+\wedge dx_1\wedge dx_2\wedge dx_3.\label{starrho1s}
\end{eqnarray}
To obtain an explicit example with non--trivial intrinsic torsion, we make the ansatz $\theta=ydz$ and $U(x,y,z)=x$ on $\{x>0\}$. The metric is therefore
$$
g=dx_1^2+dx_2^2+dx_3^2+dx_4^2+xdx^2+xdy^2+(x+\frac{y^2}{x})dz^2+\frac{1}{x}dt^2+2\frac{y}{x}dzdt
$$
with orthonormal frame
\begin{equation*}
\begin{array}{llll}
e_1=\partial_{x_1}, & e_2=\partial_{x_2}, & e_3=\partial_{x_3}, &
e_8=\partial_{x_4},\\ e_4=\frac{1}{\sqrt{x}}\partial_y, &
e_5=-\sqrt{x}\partial_t, & e_6=-\frac{1}{\sqrt{x}}\partial_x, &
e_7=\frac{1}{\sqrt{x}}(\partial_z-y\partial_t).
\end{array}
\end{equation*}
The only non--trivial brackets are
$$
\begin{array}{ll}
[e_4,e_6]=-\frac{1}{2\sqrt{x}^3}e_4, &
[e_5,e_6]=\frac{1}{2\sqrt{x}^3}e_5,\\
\lbrack e_4,e_7\rbrack=\frac{1}{\sqrt{x}^3}e_5, &
[e_6,e_7]=\frac{1}{2\sqrt{x}^3}e_7.
\end{array}
$$
Since the anti--self--dual 2--forms $\omega_1$, $\omega_2$ and
$\omega_3$ are the associated K\"ahler forms of the hyperk\"ahler
structure on $M$, we have $\nabla\omega_i=0$. From this and Koszul's formula
$$
2g(\nabla_{e_i}e_j,e_k)=g([e_i,e_j],e_k)+g([e_k,e_i],e_j)+g([e_k,e_j],e_i)=c_{ijk}+c_{kij}+c_{kji}
$$
we deduce
$$
\nabla(e_6\wedge e_7)=\nabla(e_4\wedge
e_5)=-\frac{1}{12}\cdot\frac{1}{\sqrt{x^3}}(e_4\otimes\omega_1^++e_5\otimes\omega_2^+),
$$
whence
$$
\nabla\rho = -\frac{1}{8\sqrt{3}}\cdot\frac{1}{\sqrt{x^3}}(e_4\otimes\omega_1^+\wedge
e_8+e_5\otimes\omega_2^+\wedge e_8).
$$
Note that $g$ is Ricci--flat (for $(M^4,h)$ is Ricci--flat) despite non--vanishing intrinsic torsion.

\medskip
{\bf Compact examples.} 
Consider the nilpotent Lie algebra $\mf{g}=\langle e_2,\ldots,e_8\rangle$ whose structure constants are determined by
\begin{equation}\label{structure}
de_i=\left\{\begin{array}{cl}0&\quad i=2,\ldots,7\\ e_{47}+e_{56}=\omega_1^+&\quad i=8\end{array}\right..
\end{equation}
The only non--trivial structure constants
are $c_{478}=-c_{748}=c_{568}=-c_{658}=1$.
Let $G$ be the associated simply--connected Lie group. The rationality of the structure constants guarantees the existence of a
lattice $\Gamma$ for which $N=\Gamma\backslash G$ is compact~\cite{ma51}. We let $M=T^2\times N$ with $e_i=dt_i,\,i=1,2$
on the torus, hence $de_i=0$. We take the basis $e_1,\ldots,e_8$ to
be orthonormal on $M$ and denote by $g$ the corresponding metric. As in~(\ref{rho1s}), the 3--form $\rho=e_{123}/2+\sum e_i\wedge\omega_i^-/4+\sqrt{3}e_8\wedge\omega_3^+/4$ defines a topological $PSU(3)$--structure whose invariant 5--form is given by~(\ref{starrho1s}). Then~(\ref{wedge}) and~(\ref{structure}) imply
$$
d\rho =
\frac{\sqrt{3}}{2}de_8\wedge\omega_3^+=\frac{\sqrt{3}}{2}\omega_1^+\wedge\omega_3^+=0
$$
and
\begin{eqnarray*}
d\star\rho & = & e_{4567}\wedge de_8-\frac{1}{2}\omega_{1-}\wedge
e_{23}\wedge de_8+\frac{1}{2}\omega_{2-}\wedge e_{13}\wedge
de_8\\
& &-\frac{1}{2}\omega_{3-}\wedge e_{12}\wedge de_8\\
& = & e_{4567}\wedge \omega_{1+}-\frac{1}{2}\omega_{1-}\wedge
e_{23}\wedge \omega_{1+}+\frac{1}{2}\omega_{2-}\wedge e_{13}\wedge
\omega_{1+}\\
& & -\frac{1}{2}\omega_{3-}\wedge e_{12}\wedge \omega_{1+}\\
& = & 0.
\end{eqnarray*}
Hence the $PSU(3)$--structure is harmonic. To show that the intrinsic torsion is non--trivial, we compute the covariant derivatives $\nabla e_i$, which are given by
$$
\nabla
e_i=\left\{\begin{array}{cl}0&\quad i=1,\,2,\,3\\-\frac{1}{2}(e_7\otimes
e_8+e_8\otimes e_7)&\quad i=4\\-\frac{1}{2}(e_6\otimes e_8+e_8\otimes
e_6)&\quad i=5\\\frac{1}{2}(e_5\otimes e_8+e_8\otimes
e_5)&\quad i=6\\\frac{1}{2}(e_4\otimes e_8+e_8\otimes
e_4)&\quad i=7\\\frac{1}{2}(-e_4\otimes e_7+e_7\otimes e_4-e_5\otimes
e_6+e_6\otimes e_5)&\quad i=8.\end{array}\right.
$$
Now $\nabla_{e_4}(e_8\wedge\omega_+^3)=e_{457}$; since the coefficient of $e_8\wedge\omega_3^+$ is irrational
while all the remaining ones are rational, we deduce $\nabla_{e_4}\rho\not=0$. Further, a
straightforward computation shows the diagonal of the Ricci--tensor
$Ric_{ii}=\sum_jg(\nabla_{[e_i,e_j]}e_i-[\nabla_{e_i},\nabla_{e_j}]e_i,e_j)$ to be given by
$$
Ric_{ii}=\left\{\begin{array}{rl}0&\quad i=1,\,2,\,3,\,8\\-\frac{1}{2}&\quad i=4,\,5,\,6,\,7.\end{array}\right.
$$
In particular, it follows that $(M,g)$ is of negative scalar curvature, but not Einstein, that is, $Ric$ has a non--trivial $\mb{1}$-- and $[2,2]$--component.

A non--trivial compact example of a harmonic $Sp(1)\cdot Sp(2)$--structure was given in~\cite{sa01}, where Salamon constructed a compact topological quaternionic K\"ahler 8--manifold $M$ whose structure form $\Omega$ is
closed, but not parallel. The example is of the form $M=N^6\times
T^2$, where $N^6$ is a compact nilmanifold associated with the Lie algebra given by
$$
de_i=\left\{\begin{array}{cl}0&\quad i=1,\,2,\,3,\,5\\e_{15}&\quad i=4\\e_{13}&\quad i=6.\end{array}\right.
$$
Therefore, the structure constants are trivial except for $c_{154}=-c_{514}=c_{136}=-c_{316}=1$, which implies
$$
\nabla
e_i=\left\{\begin{array}{cl}0&\quad i=2,\,4,\,6,\,7,\,8\\-\frac{1}{2}(e_3\otimes
e_6+e_4\otimes e_5+e_5\otimes
e_4+e_6\otimes e_3)&\quad i=1\\\frac{1}{2}(e_1\otimes e_6+e_6\otimes e_1)&\quad i=3\\\frac{1}{2}(e_1\otimes e_4+e_4\otimes e_1)&\quad i=5.\end{array}\right.
$$
It follows that
$$
Ric_{ii}=\left\{\begin{array}{cl}0&\quad i=2,\,4,\,6,\,7,\,8\\
-\frac{1}{2}&\quad i=1\\
-\frac{1}{4}&\quad i=3,\,5
\end{array}\right.
$$
so $(M,g)$ is of negative scalar curvature, but not Einstein, that is, $Ric$ has a non--trivial $\mathbf{1}$-- and $[1,1,2]$--component. Summarising, we obtain the

\begin{thm}\label{riccob}
\hspace{-2pt}Compact harmonic $PSU(3)$-- and $Sp(1)\cdot Sp(2)$--structures with non--vanishing intrinsic torsion do exist. Further, they are not necessarily Einstein. 
\end{thm}
%
%
%
%
%
\appendix\section{A matrix representation of $\cliff(\R^8,g_0)$}
\label{matrix}
For a fixed orthonormal basis $e_1,\ldots, e_8$ of $(\Lambda^1,g)\cong(\R^8,g_0)$ let $E_{ij}=e_i\wedge e_j$ denote the basis of $\Lambda^2$ which we identify with skew--symmetric matrices via
\begin{eqnarray*}
E_{ij} & = & \left(\begin{array}{ccccc} 0 & \ldots & \ldots & \ldots & 0\\
\ldots & \ldots & \ldots & -1 & \ldots\\ \ldots & \ldots & \ldots & \ldots & \ldots\\
\ldots & 1 & \ldots & \ldots & \ldots \\ 0 & \ldots & \ldots &
\ldots & 0\end{array}\right)
\begin{array}{c} \\ \ldots i \\ \\ \ldots j\\ \\ \end{array}.\\[-5pt]
& &\qquad\qquad\mspace{-10pt}\vdots \qquad \quad\mspace{1pt}\vdots\\[-7pt]
& &\qquad\qquad\mspace{-10pt} i \qquad \quad \mspace{1pt} j
\end{eqnarray*}
Then the matrix representation $\kappa:\cliff(\R^8,g_0)\to\End(\Delta_+\oplus\Delta_-)$ computed from~(\ref{cliff(8)}) with respect to the standard basis $e_1=1,\,e_2=i,\,e_3=j,\,e_4=k,\,e_5=e,\,e_6=e\cdot i,\,e_7=e\cdot j,\,e_8=e\cdot k$
of $(\O,\enorm{\cdot})$ is given by 
\begin{equation}\label{spin8rep}
\begin{array}{lcr}
\kappa(e_1) & = & -E_{1,9}-E_{2,10}-E_{3,11}-E_{4,12}-E_{5,13}-E_{6,14}-E_{7,15}-E_{8,16},\\
\kappa(e_2) & = & \phantom{-}E_{1,10}-E_{2,9}-E_{3,12}+E_{4,11}-E_{5,14}+E_{6,13}+E_{7,16}-E_{8,15},\\
\kappa(e_3) & = & \phantom{-}E_{1,11}+E_{2,12}-E_{3,9}-E_{4,10}-E_{5,15}-E_{6,16}+E_{7,13}+E_{8,14},\\
\kappa(e_4) & = & \phantom{-}E_{1,12}-E_{2,11}+E_{3,10}-E_{4,9}-E_{5,16}+E_{6,15}-E_{7,14}+E_{8,13},\\
\kappa(e_5) & = & \phantom{-}E_{1,13}+E_{2,14}+E_{3,15}+E_{4,16}-E_{5,9}-E_{6,10}-E_{7,11}-E_{8,12},\\
\kappa(e_6) & = & \phantom{-}E_{1,14}-E_{2,13}+E_{3,16}-E_{4,15}+E_{5,10}-E_{6,9}+E_{7,12}-E_{8,11},\\
\kappa(e_7) & = & \phantom{-}E_{1,15}-E_{2,16}-E_{3,13}+E_{4,14}+E_{5,11}-E_{6,12}-E_{7,9}+E_{8,10},\\
\kappa(e_8) & = &
\phantom{-}E_{1,16}+E_{2,15}-E_{3,14}-E_{4,13}+E_{5,12}+E_{6,11}-E_{7,10}-E_{8,9}.
\end{array}
\end{equation}
%
%
%
%
%
\section{A matrix representation of the invariant supersymmetric maps}
\label{susymatrix}
1. $PSU(3)$
\newline
With respect to a $PSU(3)$--frame and a fixed orthonormal basis of $\Delta_{\pm}$, the invariant supersymmetric maps $\sigma_{\pm}:\Lambda^1\to\Delta_{\pm}$ are given (up to a scalar) by
$$
\sigma_+ = \left(\begin{array}{rrrrrrrr}
\scriptstyle0 &\scriptstyle -\frac{1}{2} &\scriptstyle 0 &\scriptstyle \frac{1}{4} &\scriptstyle -\frac{\sqrt{3}}{4} &\scriptstyle -\frac{1}{4} &\scriptstyle \frac{\sqrt{3}}{4} &\scriptstyle \frac{1}{2}\\[3pt]

\scriptstyle-\frac{1}{2} &\scriptstyle 0 &\scriptstyle -\frac{1}{2} &\scriptstyle -\frac{\sqrt{3}}{4} &\scriptstyle -\frac{1}{4} &\scriptstyle -\frac{\sqrt{3}}{4} &\scriptstyle -\frac{1}{4} &\scriptstyle 0\\[3pt]

\scriptstyle0 &\scriptstyle \frac{1}{2} &\scriptstyle 0 &\scriptstyle \frac{1}{4} &\scriptstyle -\frac{\sqrt{3}}{4} &\scriptstyle \frac{1}{4} &\scriptstyle -\frac{\sqrt{3}}{4} &\scriptstyle \frac{1}{2}\\[3pt]

\scriptstyle-\frac{1}{2} &\scriptstyle 0 &\scriptstyle \frac{1}{2} &\scriptstyle \frac{\sqrt{3}}{4} &\scriptstyle \frac{1}{4} &\scriptstyle -\frac{\sqrt{3}}{4} &\scriptstyle -\frac{1}{4} &\scriptstyle 0\\[3pt]

\scriptstyle0 &\scriptstyle -\frac{1}{2} &\scriptstyle 0 &\scriptstyle -\frac{1}{4} &\scriptstyle \frac{\sqrt{3}}{4} &\scriptstyle \frac{1}{4} &\scriptstyle -\frac{\sqrt{3}}{4} &\scriptstyle \frac{1}{2}\\[3pt]

\scriptstyle\frac{1}{2} &\scriptstyle 0 &\scriptstyle \frac{1}{2} &\scriptstyle -\frac{\sqrt{3}}{4} &\scriptstyle -\frac{1}{4} &\scriptstyle -\frac{\sqrt{3}}{4} &\scriptstyle -\frac{1}{4} &\scriptstyle 0\\[3pt]

\scriptstyle0 &\scriptstyle -\frac{1}{2} &\scriptstyle 0 &\scriptstyle \frac{1}{4} &\scriptstyle -\frac{\sqrt{3}}{4} &\scriptstyle \frac{1}{4} &\scriptstyle -\frac{\sqrt{3}}{4} &\scriptstyle -\frac{1}{2}\\[3pt] 

\scriptstyle-\frac{1}{2} &\scriptstyle 0 &\scriptstyle \frac{1}{2} &\scriptstyle -\frac{\sqrt{3}}{4} &\scriptstyle -\frac{1}{4} &\scriptstyle \frac{\sqrt{3}}{4} &\scriptstyle \frac{1}{4} &\scriptstyle 0
\end{array}\right)\, \sigma_- = \left(\begin{array}{rrrrrrrr}
\scriptstyle-\frac{1}{2} &\scriptstyle 0 &\scriptstyle \frac{1}{2} &\scriptstyle \frac{\sqrt{3}}{4} &\scriptstyle -\frac{1}{4} &\scriptstyle -\frac{\sqrt{3}}{4} &\scriptstyle \frac{1}{4} &\scriptstyle 0\\[3pt]

\scriptstyle0 &\scriptstyle -\frac{1}{2} &\scriptstyle 0 &\scriptstyle \frac{1}{4} &\scriptstyle \frac{\sqrt{3}}{4} &\scriptstyle \frac{1}{4} &\scriptstyle \frac{\sqrt{3}}{4} &\scriptstyle \frac{1}{2}\\[3pt]

\scriptstyle-\frac{1}{2} &\scriptstyle 0 &\scriptstyle -\frac{1}{2} &\scriptstyle -\frac{\sqrt{3}}{4} &\scriptstyle \frac{1}{4} &\scriptstyle -\frac{\sqrt{3}}{4} &\scriptstyle \frac{1}{4} &\scriptstyle 0\\[3pt]

\scriptstyle0 &\scriptstyle \frac{1}{2} &\scriptstyle 0 &\scriptstyle \frac{1}{4} &\scriptstyle \frac{\sqrt{3}}{4} &\scriptstyle -\frac{1}{4} &\scriptstyle -\frac{\sqrt{3}}{4} &\scriptstyle \frac{1}{2}\\[3pt]

\scriptstyle-\frac{1}{2} &\scriptstyle 0 &\scriptstyle \frac{1}{2} &\scriptstyle -\frac{\sqrt{3}}{4} &\scriptstyle \frac{1}{4} &\scriptstyle \frac{\sqrt{3}}{4} &\scriptstyle -\frac{1}{4} &\scriptstyle 0\\[3pt]

\scriptstyle0 &\scriptstyle \frac{1}{2} &\scriptstyle 0 &\scriptstyle \frac{1}{4} &\scriptstyle \frac{\sqrt{3}}{4} &\scriptstyle \frac{1}{4} &\scriptstyle \frac{\sqrt{3}}{4} &\scriptstyle -\frac{1}{2}\\[3pt]

\scriptstyle\frac{1}{2} &\scriptstyle 0 &\scriptstyle \frac{1}{2} &\scriptstyle -\frac{\sqrt{3}}{4} &\scriptstyle \frac{1}{4} &\scriptstyle -\frac{\sqrt{3}}{4} &\scriptstyle \frac{1}{4} &\scriptstyle 0\\[3pt]

\scriptstyle0 &\scriptstyle \frac{1}{2} &\scriptstyle 0 &\scriptstyle-\frac{1}{4} &\scriptstyle -\frac{\sqrt{3}}{4} &\scriptstyle \frac{1}{4} &\scriptstyle \frac{\sqrt{3}}{4} &\scriptstyle \frac{1}{2} 
\end{array}\right)
$$

2. $Sp(1)\cdot Sp(2)$
\newline
With respect to an $Sp(1)\cdot Sp(2)$--frame and a fixed orthonormal basis of $\Delta_+$, the invariant supersymmetric map $\sigma_+:\Lambda^1\to\Delta_+$ is given (up to a scalar) by
$$
\sigma_+=\left(
\begin {array}{rrrrrrrr} \scriptstyle\phantom{-}1&\scriptstyle0&\scriptstyle0&\scriptstyle0&\scriptstyle0&\scriptstyle0&\scriptstyle0&\scriptstyle0
\\[3pt]
\scriptstyle0&\scriptstyle
-1&\scriptstyle0&\scriptstyle0&\scriptstyle0&\scriptstyle0&\scriptstyle0&\scriptstyle0
\\[3pt]
\scriptstyle0&\scriptstyle0&\scriptstyle-1&\scriptstyle0&\scriptstyle0&\scriptstyle0&\scriptstyle0&\scriptstyle0
\\[3pt]
\scriptstyle0
&\scriptstyle0&\scriptstyle0&\scriptstyle-1&\scriptstyle0&\scriptstyle0&\scriptstyle0&\scriptstyle0
\\[3pt]
\scriptstyle0&\scriptstyle0&\scriptstyle0&\scriptstyle0&\scriptstyle\phantom{-}1&\scriptstyle0&\scriptstyle0&\scriptstyle0
\\[3pt]
\scriptstyle0
&\scriptstyle0&\scriptstyle0&\scriptstyle0&\scriptstyle0&\scriptstyle\phantom{-}1&\scriptstyle0&\scriptstyle0
\\[3pt]
\scriptstyle0&\scriptstyle0&\scriptstyle0&\scriptstyle0&\scriptstyle0&\scriptstyle0&\scriptstyle\phantom{-}1&\scriptstyle0
\\[3pt]
\scriptstyle0
&\scriptstyle0&\scriptstyle0&\scriptstyle0&\scriptstyle0&\scriptstyle0&\scriptstyle0&\scriptstyle\phantom{-}1
\end{array}\right)
$$
{\bf Acknowledgments}

This paper is based on a part of the author's doctoral thesis~\cite{wi05}. He wishes to acknowledge the EPSRC, the University of Oxford, the DAAD and the Studienstiftung des deutschen Volkes for various forms of financial support. Further, he would like to thank the \'Ecole Polytechnique, Palaiseau for a research assistantship whilst the preparation of this paper. The author is also happy to acknowledge Michael Crabb for very valuable ideas and the thesis' examiners Simon Salamon and Dominic Joyce for their constructive comments. Finally, he wishes to express his special gratitude to his supervisor Nigel Hitchin.

\end{document}